\documentclass{amsart}
\usepackage{amsmath, amsthm}
\usepackage{amssymb, latexsym}
\usepackage{amsfonts}

 \usepackage{graphicx}
\usepackage[dvipsnames]{xcolor}

\usepackage{tkz-graph}

\usepackage{calrsfs}

\newtheorem{theorem}[equation]{Theorem}

\newtheorem{proposition}[equation]{Proposition}
\newtheorem{lemma}[equation]{Lemma}

\newtheorem{corollary}[equation]{Corollary}

\newcounter{com}

\font\Bbb=msbm10
\def\semi{\hbox{\Bbb o}}

\def\Z{\hbox{\Bbb Z}}

\numberwithin{equation}{section}

\let\define=\def
\let\redefine=\def

\def\Aut{{\rm Aut}}

\def\rank{{\rm rank}}

\define\a{{ \alpha }}
\redefine\b{{ \beta }}
\redefine\c{{ \gamma }}
\redefine\d{ \delta }
\define\e{{ \epsilon }}
\redefine\g{\gamma}

\redefine\D{{ \Delta }}

\define\z{{ \zeta }}

\define\({ \left( }
\define\){ \right) }
\define\[{ \left[ }
\define\]{ \right] }
\define\<{ \langle }
\define\>{ \rangle }
\let\ljunk=\{
\let\rjunk=\}
\redefine\{{\left\ljunk}
\redefine\}{\right\rjunk}

\redefine\.{ \diamomd }
\redefine\-{ \ominus }
\redefine\+{\oplus}
\redefine\.{ \cdot}
\redefine\#{\sharp}

\redefine\.{ \cdot }

        \def\rad{{\rm rad}}

        \def\GF{{\rm GF}}

     \def\diag{{\rm diag}}

        \font\Aaa=msam10
        \def\nor{\hbox{\Aaa C}}
        \def\ron{\hbox{\Aaa B}}
        \def\noreq{\mathop{\hbox{\Aaa E}}}
        \def\roneq{\hbox{\Aaa D}}

\def\half{\frac{1}{2}}

        \def\wt{{\rm wt}}
\let\Gamma\varGamma
\let\Delta\varDelta
\let\Theta\varTheta
\let\Lambda\varLambda
\let\Xi\varXi
\let\Pi\varPi
\let\Sigma\varSigma
\let\Upsilon\varUpsilon
\let\Phi\varPhi
\let\phi\varphi
\let\Psi\varPsi
\let\Omega\varOmega

     \def\ell{{l}}

\def\ndiv{ {\not\kern-.5pt\hbox{\big |}\,} }

\def\Syl{{\rm Syl}}

\begin{document}

\title[Collineation Groups]
{Antiflag Transitive Collineation Groups }

       \author{William  M. Kantor
      }
       \address{U. of Oregon,
       Eugene, OR 97403
        \      and  \
       Northeastern U., Boston, MA 02115}
       \email{kantor@uoregon.edu}

\begin{abstract}  
We present   known results concerning antiflag 
transitive collineation groups of finite
projective spaces and finite polar spaces.

\end{abstract}

\maketitle

 \vspace{-10pt}
 
 \section*{Preface}
This is  a revision of a paper 
by P. J. Cameron and W. M. Kantor, 
``$2$-Transitive and Antiflag Transitive Collineation Groups of Finite
Projective Spaces'',  {J. Algebra}  {60} (1979)  384--422.
All theorems in that paper are corrected by adding further   imprimitive antiflag transitive subgroups  
 to various statements:    
 \begin{itemize}
\item 
$ SL(\half n,16)  $ or   $Sp(\half n,16)   \nor G <\Gamma L(n,4)$,
 or
 $ G_2(16)  \nor  G < \Gamma L(12,4)$;\vspace{2pt}
\item 
$Sp(\half n,16)   \nor G <\Gamma Sp(n,4) \cong \Gamma O(n+1,4)$,
 or\vspace{2pt}
 $ G_2(16)  \nor  G < \Gamma  Sp(12,4) \cong\vspace{2pt} \Gamma O(13,4)$; 
\item 
$SU(\half n, 4) \nor  G <\Gamma O^\epsilon(n, 4)$, where $\epsilon = (-)^{\half n};$
 and 
\item 
$G\cong A_9$  inside   $\Omega^+(8,2).$ 
\end{itemize}

  This revision  uses the same methodology as the original paper.  In particular, 
  it does  not involve  more recent group theory.   

  The actual results can be deduced from   more recent results that depend on the 
Classification of the Finite Simple Groups.   Liebeck \cite{Liebeck}  completed results of   Hering concerning transitive finite linear groups,
so that paper implies Theorems I-III.  The transitivity results in Theorems IV-V are implicitly very special cases of  Guralnick-Penttila-Praeger-Saxl \cite{GPPS}.

However, 
I strongly believe that the aforementioned Classification should not~be invoked when it is not needed.
Moreover,   there are    surprising and entertaining~parts of the original proofs
(especially the appearance of generalized hexagons).   
 The  notation, methodology  and  relatively  minimal background of the original paper are  retained here.
   As much as possible   
    the original paper has been left unchanged; for example, the  numbering  of intermediate results is not altered.
   
\section{Introduction}
An unpublished result of Perin [20]  states that a subgroup of $\Gamma L(n, q), n\ge 3,$
that induces a primitive rank 3 group of even order on the set of points of
$PG(n - 1, q)$, necessarily preserves a symplectic polarity. (Such groups are
essentially known, if $q>3$, by another theorem of Perin   [19].) The present paper extends
both Perin's result and his method, in order to deal with some familiar problems
concerning collineation groups of finite projective spaces; among these, 2-transitive collineation groups  [25], and   both the case of  semilinear groups and
 the case $q \le 3 $ of  Perin's theorem [19].

An {\em antiflag} is an ordered pair consisting of a hyperplane and a point not on it;
if the underlying vector space is endowed with a symplectic, unitary or orthogonal
geometry, both the point and the pole of the hyperplane are assumed to be
isotropic or singular. Our main results are the following four theorems.

{
\renewcommand{\theequation}{I}
\begin{theorem}
If $G \le \Gamma L(n, q), $  $n\ge 3,$ and $G$ is $2$-transitive on the set of
points of $PG(n - 1, q),$  then either $G \ge SL(n, q),$  or $G$ is $A_7$ inside $SL(4, 2).$

\end{theorem}
}

{
\renewcommand{\theequation}{II}
\begin{theorem}
If $G \le \Gamma L(n, q)$  and $G$ is transitive on antiflags and primitive
but not $2$-transitive on points$,$  then $G$ preserves a symplectic polarity$,$ 
 and one of the  following holds$:$
\begin{itemize}
\item[\rm (i)]  
$ G \roneq  Sp(n, q);$
\item[\rm (ii)]$G$ is $A_6$ inside $Sp(4, 2); $  or
\item[\rm (iii)] $G \roneq G_2(q), $  $q$ even$,$ and $G$ acts on the generalized hexagon associated with $ G_2(q),$  which is itself embedded naturally in $PG(5, q).$
\end{itemize}
\end{theorem}
}

{
\renewcommand{\theequation}{III}
\begin{theorem}
 If  $G \le \Gamma L(n, q)$ and $G$ is transitive on antiflags and imprimitive
on points$,$ then 
$q = 2$ or $4 $  and $ \Gamma L(\frac{1}{2}n, q^2)\ge G \ron SL(\frac{1}{2}n, q^2), $ $Sp(\frac{1}{2}n, q^2),$ or  $ G_2(q^2)$
$($with $n = 12).$%
\footnote{\,We are grateful to  Nick Inglis and Jan Saxl   for pointing out that the 
case  $  \Gamma L(\half n, 16) $ had been omitted.  This  led to the other groups
over $GF(16)$.}%
\   In each case$,$  $G$ is embedded naturally in $\Gamma L(n, q).$ 
 \end{theorem}
}

{
\renewcommand{\theequation}{IV}
\begin{theorem}
  If $G \le \Gamma Sp(n, q), $ $ \Gamma O^\pm(n, q)$ or  $ \Gamma U(n, q),$ for a classical
geometry of rank at least $3,$ and $G$ is transitive on antiflags$,$ then one  of the following
holds $($and the embedding of $G$ is the natural one$)$$:$
  
  \begin{itemize}
\item[\rm (i)]  
$G \roneq Sp(n, q),$  $ \Omega^\pm(n, q),$ resp. $ S U(n, q);$
 
\item[\rm (ii)] $G \roneq G_2(q)$ inside $\Gamma O(7, q) $  
$($or $\Gamma Sp(6, q), $  $q$ even$);$

\item[\rm (iii)]  $\Omega (7, q) \noreq   G/Z(G) \!< \! P \Gamma O^+(8, q), \!$ with $ G/Z(G)$ conjugate in  $\Aut(P\Omega^+(8, q) ) $ to a group fixing a nonsingular $1$-space$;$

\item[\rm (iv)]   For $q=2$ or $4 ,$ $Sp(\half n, q^2) \nor G < \Gamma Sp(n, q)\cong \Gamma O(n+1,q);$$\,{}^2$

\item[\rm (v)]   For $q=2$ or $4,$
$G_2(q^2)\nor G  <  \Gamma  Sp(12, q) \cong \Gamma O(13, q);$$\,{}^2$
 \item[\rm (vi)]  
   For $q=2$ or $4,$
$SU(\half n, q) \nor  G <\Gamma O^\epsilon(n, q),$ where $\epsilon = (-)^{\half n};\,$\footnote{\,The case $q=4$ had been omitted.}
or
\item[\rm (vii)]
$G\cong A_9 $  inside $\Omega^+(8,2)$.%
\footnote{\,This case had been omitted.}

\end{itemize} 

\end{theorem}
}

Theorem I solves a problem posed by Hall and Wagner [25], which  has been
studied by Higman [8, 10], Perin [19], Kantor [13] and Kornya [15].
An  independent and alternative approach to this theorem is given by Orchel [16]; we are
grateful to Orchel for sending us a copy of his thesis.

If $G$ is 2-transitive, then $G$ is antiflag transitive; and also $G_H^H$ is  
antiflag transitive for each hyperplane $H$. This elementary fact allows us to use induction.
(Indeed, Theorems I-III are proved simultaneously by induction in Part I of
this paper.) 
 The groups in Theorem III and Theorem IV(iv-vi) must contain both the indicated
 quasisimple group and $\Aut(GF(q^2))$.
 Another problem, solved in Theorems II and IV, is that of primitive
rank 3 subgroups of classical groups. This was posed by Higman and
McLaughlin [11], and solved by Perin [19] 
(for linear groups) and Kantor and Liebler [14] except
in the cases $Sp(2 m, 2)\cong \Omega(2 m+ 1, 2)$ and $Sp(2 m, 3)$. Here, induction is made possible
by  fact that the stabilizer of a point $x$ is antiflag transitive   on  $x^\perp/x$.  

The striking occurrence of  $G_2(q)$ in these theorems is related to  a crucial
element of our approach. This case is obtained from a general embedding
theorem for metrically regular graphs (3.1), in which the  Feit-Higman theorem
[7] on generalized polygons arises unexpectedly but naturally. Other familiar
geometric objects and theorems come into play later on: the characterizations 
of projective spaces due to Veblen and Young [24] and Ostrom and Wagner [18],
as well as translation planes, arise in Theorem III, while Tits' classification of
polar spaces [23] and the triality automorphism of $P\Omega^+(8, q)$ are used for
Theorem IV.

All  of  the proofs require familiarity with the geometry of the classical groups.
On the other hand, group-theoretic classification theorems have been 
avoided. Moreover, knowledge of $G_2(q)$  is not assumed for Theorem I, and what
is required for Theorems II-IV is contained in the Appendix, where we have
given a new and elementary proof of the existence of the generalized hexagons
of type  $G_2(q).$

This paper began as an attempt to extend Perin's result [20] to rank 4
subgroups of classical groups. As in Perin [19], one case with $q = 2$  is left open:

{
\renewcommand{\theequation}{V}
\begin{theorem}
Suppose $G \le \Gamma Sp(n, q) ~(n\ge 6), $ $\Gamma O^\pm(n, q)~(n\ge7),$  or $\Gamma U(n, q)$
$(n\ge 6)$. If $G$ induces a primitive rank $4$ group on the set of isotropic or singular
points$,$  then one of the following holds$:$
\begin{itemize}
\item[\rm (i)]  
 $G\roneq  G_2(q)$ is embedded naturally in $ \Gamma O(7, q)$  $ ($or $\Gamma Sp(6, q),$
 $  q$ even$);$
\item[\rm (ii)] 
 $G \roneq \Omega(7, q), $  $q$ even$,$ or $2. \Omega(7, q), $ $q$ odd$,$   
 each  embedded irreducibly in\break
  $\Gamma O^+(8, q);$ or
\item[\rm (iii)]$G\le O^\pm(2m, 2),$ and G is transitive on the pairs $(x, L)$  with $L$ a totally singular line and $x$ a point of  $L.$
\end{itemize}
\end{theorem}
}

The examples (ii) (and (iii) in Theorem IV) are obtained by applying the
triality automorphism to the more natural $\Omega (7, q)$ inside $P \Omega^+(8, q).$  
 As for (iii), examples are $A_7$ and $S_7$ inside $O ^+(6, 2).$

Other results in a similar spirit are given in Section 8, as corollaries to Theorem~I.

Some further results are of interest independent of their application to the
above theorems. A general result on embedding metrically regular graphs in
projective spaces is proved in Section 3; this is used several times, and is crucial for all of  the theorems.
Theorem 10.3 characterizes nonsingular quadrics of dimension $2 m - 1$ contained
in an $O^+(2 m , q)$ quadric for $m\ge 3$. In Section 12, parameter restrictions are
obtained for rank 4 subgroups of rank 3 groups (and their combinatorial analogues).
Finally, the Appendix gives an elementary construction and characterization
of the $G_2(q)$ hexagon.

The paper falls into two parts. The first (Sections 2-8) deals with antiflag
transitive collineation groups of projective spaces (Theorems I-III); we note
that Sections 3 and 5, on the primitive, not 2-transitive case, are virtually self-contained.
The second part (Sections 9-14) contains the proofs of Theorems IV
and V, concerning polar spaces.

\medskip
 
\centerline{I.  THEOREMS I-III}
 
\smallskip

\section{Preliminaries}

A {\em  point} ({\em hyperplane}) of a vector space $V$ is a subspace of dimension
1 (codimension 1). If 
$V$ is $n$-dimensional over $GF(q)$, the set of  points (equipped with the structure of projective geometry)  is denoted by $PG(n - 1, q)$; but in this paper,
its dimension will always be $n$. The notation
$  SL(V) = SL(n, q), $ $GL(n, q) $ and $\Gamma  L(n, q)$ is standard.

If, in addition,  $V$ is equipped with a symplectic, unitary or orthogonal  geometry, then 
$\Gamma Sp(n, q)$, $\Gamma U(n, q)$ and   $\Gamma O^\pm (n, q)$ denote the groups of semilinear maps preserving the geometry projectively. For example. $\Gamma O^\pm(n, q)$
consists  of all invertible semilinear maps $g$ such that  $\phi(v^g) = c\phi (v)^\sigma$ for all $v\in V  $, where
$\phi$ is the quadratic form defining the geometry, $c$ is a scalar, and $\sigma $ 
is a field automorphism.
The groups $Sp(n, q),$ $ SU(n, q)$ and $\Omega^\pm(n, q)$ are defined as   usual.  
We use   totally isotropic or totally singular (abbreviated
t.i.~or t.s.) subspaces of these  geometries.
There is some ambiguity in the terminology ``t.i.~or t.s.~subspace''   since orthogonal geometries
have both types of subspaces in characteristic 2; but in this case we aways refer to t.s.~subspaces.
We
will occasionally require the fact that $Sp(2n, q)\cong \Omega(2n +1, q)$  when $ q$ is even.
(Explicitly, if $V$ is the natural $Sp(2n, q)$-module,  then there is a 
nondegenerate $2n +1$-dimensional orthogonal space $\tilde V$ such that $\tilde V/\rad \tilde V = V$, with the natural 
map  $ \tilde V\to V$
inducing a bijection between singular and isotropic points.) 
The reader is
referred to Dieudonn\'e [6] for further information concerning all of these groups.

Points will be denoted $x, y, z$, lines $ L, L'$ and hyperplanes $H, H'$. We will 
generally identify a subspace $\Delta$  of $V$ with its set of points;  $|\Delta|$ denotes 
its number of points, and $x\in \Delta$  will be used instead of 
$x\subseteq \Delta$. Similarly, for subspaces  $ \Delta$  and $\Sigma$,
$ \Delta-\Sigma$ denotes the set of points in $ \Delta$ but not $\Sigma$.
On the other hand, 
 the dimension $\dim \Delta$ of a subspace denotes the vector space dimension.
If $A\le G L(V)$ and  $W$ is a subspace of $V$ then 
 $C_W(A)=\{w\in W\mid w^a=w,\, a\in A\}$ and 
 $[W,A]=\{w^a- w \mid w\in W,\, a\in A\}$
 are vector subspaces   that will be studied as sets of points; 
we expect that the context will make it clear whether 
a subspace  is being viewed as a set of points after being obtained as  a set of vectors.
 
We generally consider  semilinear groups;  but when discussing transitivity we always consider
 the induced (projective) group on 1-spaces (points) rather than transitivity on vectors.
If $ \Delta$  is any subset of $V$, then $G_ \Delta$  and $C_G( \Delta )$  are 
respectively the setwise and
vector-wise stabilizers of  $\Delta$  in the semilinear  group $G$;
$ G _ {\Delta \Sigma} =  G _ {\Delta} \cap  G _ {\Sigma} $.
Moreover, $G_\Delta^\Delta$ is the semilinear group   induced on $\Delta$ if  $\Delta$ is a
subspace; this group will usually be viewed projectively.
 Similarly, if $x \in H$, then $G_{xH}^{H/x}$  is the group induced  by $G_{xH} $ 
on  the space  $H/x$.

The {\em rank} of a transitive permutation group is the total number of orbits of 
the stabilizer of a point.

The remainder of this section lists further definitions and results required in
the proofs of Theorems I-V.

\begin{theorem} [{Ostrom-Wagner [18], Ostrom [17]}]
If a projective  plane  $P$ of prime power order $q$ admits a collineation group $G$   transitive on non-incident  point-line pairs$,$  then $P$  is desarguesian and $G\ge PSL(3, q).$
\end{theorem}

Of course, (2.1) is true  without the prime power assumption, but we
only need the stated case, which is much easier to prove.
 The next result is needed for (2.1), and is also used elsewhere in our argument.

\begin{theorem} [{[4, pp. 122, 130-34]}]
Let ${\mathfrak  A}$ be an affine translation  plane  of order $q,$\, $L$ a line$,$ $x\in L,$ and $E$ the group of elations  with center $x$ and axis $L$. Then 
\begin{itemize}
\item[\rm (i)] $E$ is semiregular on the set of lines different from $L$ on $x;$
and
 \item[\rm (ii)]  If $|E|=q$ for each $L$ and $x,$ then ${\mathfrak  A}$  is desarguesian.
\end{itemize}
\end{theorem}

Additional, more elementary results concerning translation planes will also be
required; the reader is referred to Dembowski [4, Chap. 4] for further information
concerning perspectivities and Baer involutions.

Consider next a geometry ${\mathcal G}$  of points, with certain subsets called ``lines",
such that any two points are on at most one line, each line has at least three
points, and each point is on at least three lines. Call ${\mathcal P}$  and ${\mathcal L}$
the sets of points
and lines. If  $a, b\in {\mathcal P}\cup {\mathcal L} ,$  the distance $\partial (a, b)$ between them is the smallest
number $k$ for which there is a sequence $a = a_0 , a_1 , \dots, a_k = b$, with each
$a_i \in {\mathcal P}\cup {\mathcal L} $   and $a_i$ incident with $a_{i+1}$  for $ i = 0, ... , k - 1.$ Such a sequence is called a ``path" from $a$ to $b$. Now ${\mathcal G}$  is a 
{\em generalized $n$-gon} ($n \ge 3$) if
\begin{itemize}
\item[\rm (i)]  whenever $\partial(a,b)<n$ there is a unique shortest path from $a$ to $b$;
 \item[\rm (ii)] for all $a$ and $ b,$  $\partial(a,b)\le n;$
and
  \item[\rm (ii)] there exist  $a$ and $ b $ with   $\partial(a,b)= n.$
\end{itemize}

A generalized $n$-gon has {\em parameters} $s,t$ if each line has exactly $s+1$ points and each point is on exactly $t+1$ lines.

\begin{theorem} [Feit-Higman {[7]}]
Generalized $n$-gons can exist only for $n=3,4,6$ or $8;$
those with $n=8$ cannot have parameters $s,s$.
\end{theorem}

Generalized quadrangles enter our considerations as the geometries of points
and lines in low-dimensional symplectic, unitary, and orthogonal geometries.
Generalized hexagons are much less familiar; the ones we need are discussed in
the Appendix (see also Sections 3,  5 below).

Generalized $n$-gons are special cases of {\em metrically regular} graphs.
Let $\Gamma$ be a connected graph defined on a set $X$ of vertices.
If $x,y\in X$, let $d(x,y)$ denote the distance between them.  Let $d$ be the diameter, and   
 $\Gamma_i(x) $ the set of points at distance $i$ from $x$, for $0\le i\le d$.  Then 
  $\Gamma$ is  metrically regular  if
 
\begin{itemize}
\item[\rm (i)]   $| \Gamma_i(x) |$ depends only on $i,$ not on $x;$ and  
 \item[\rm (ii)] if $d(x,y)=i,$ the numbers of points at distance $1$ from $x$ and distance 
 $i-1 ~($resp. $i,$ $i+1)$ from $y$ depend only on $i,$ and not on $x$ and $y$.
\end{itemize}
\smallskip

(Condition  (i)  follows from  (ii)  here.)

If  ${\mathcal G}$  is a geometry 
as previously defined, its {\em point graph}  $\Gamma$  is obtained by
joining two points of   ${\mathcal G}$  by an edge precisely when they are distinct and collinear.
This graph may be metrically regular; for example, it is so when  ${\mathcal G}$   is a generalized
$n$-gon. (Here the distances $d$ and $\partial $  in graph and geometry are related by 
$d(x, y)= \frac{1}{2}\partial(x, y)$ for $x, y \in{\mathcal P}.)$

If $n$ is an integer then $n_p$ denotes the largest power of $p$ dividing $n$
 (where $p$, as always, is a prime).

If $q$ is a power $p^e$ of $p$, and $ k\ge 2$, a 
{\em primitive divisor}
 of   $q^k-1$ is a prime $r\mid q^k-1$ such that $r\not| \,p^i-1$ for $1<p^i<q^k$.
Note that   $r\equiv  1	$ (mod $ek$), by Fermat's theorem.

\begin{theorem} [Zsigmondy {[28]}]
 If $q>1$ is a power of $p$ and $k>1,$ then $q^k-1$ has a primitive divisor unless either 
 
\begin{itemize}
\item[\rm (i)]    $k = 2$ and $q$ is a Mersenne prime$,$ or 
 \item[\rm (ii)] $ q^k = 64.$
 \end{itemize}
 
\end{theorem}

\section{Embedding Metrically regular graphs in projective spaces}

In this section we will prove a general result concerning certain embeddings
in projective spaces. Let ${\mathcal G}$ be a geometry, with point set $\Omega$  and point graph $\Gamma$.
For $x \in \Omega$, let $W_i(x)$  be the set of points distant at most $i$ from $x$.  We assume
the following axioms (for all $x \in \Omega$):
\begin{itemize}
\item[\rm (a)]    $\Omega$ is a set of points spanning  $PG(n -1,q);$
 \item[\rm (b)] each line $L$ of ${\mathcal G}$ (or ${\mathcal G}$-{\em line}) is a line of $PG(n - 1,q);$

  \item[\rm (c)]$\Omega$ is the union of the set of  ${\mathcal G}$-lines;
   \item[\rm (d)]$\Gamma$  is metrically regular with diameter $d\ge 2;$
    \item[\rm (e)]$W_1(x)$  is a subspace of $PG(n - 1, q);$
     \item[\rm (f)]$W_i(x) = \Omega\cap  U_i(x)$ for some subspace $ U_i(x)$
       for each $i;$ and
      \item[\rm (g)]$|W_2(x)| = (q^h - 1  )/(q - 1)$  for some  $h$.
     \end{itemize}

Note that (a)-(d) are among the embedding hypotheses in Buekenhout-Lef\`evre [1].

In (3.1) and (3.2)  we will determine all geometries satisfying (a)-(g).  For
Theorem I, a complete classification is not required; the weaker result  (3.1) suffices.

\begin{theorem} 
If ${\mathcal G}$ satisfies  {\rm(a)-(g),}  then either 
\begin{itemize}
\item[\rm (i)]    $d = 2$ and ${\mathcal G}$ consists of the totally isotropic points and lines of a symplectic
polarity $x\leftrightarrow W_1(x);$   or
 \item[\rm (ii)] 
 $d = 3,$   ${\mathcal G}$ is a generalized hexagon with parameters $q, q,$
 and    each  $W_1(x)$ has dimension $3$.
 $($Moreover$,$ 
if $W_2(x)$ and $W_3(x)$ are subspaces for all $x,$ then $n = 6$ and 
$x\leftrightarrow W_2(x) $   is a symplectic polarity.$)$
      \end{itemize}
\end{theorem}
\proof
Set $m = \dim W_1(x)$ (recalling from Section 2 that   ``dim" means
vector space dimension).  If    $d(x,y)=i \ge 1$, let
$$
\begin{array}{lllll}
e_i\hspace{-6pt}&=  \dim W_1(x) \cap W_{i-1} (y) ,  \smallskip
\\
f_i\hspace{-6pt}&=\dim W_1(x) \cap W_{i}  (y) .
\end{array}
$$
(Note that both $W_1(x) \cap W_{i-1} (y) $ and $W_1(x) \cap W_{i} (y) $ are subspaces.  For if $W_j(y)=\Omega\cap U_j(y)$, then 
$W_1(x) \cap W_{j}  (y)   =   W_1(x) \cap \Omega\cap U_j(y) = W_1(x) \cap U_j(y).$)
These dimensions depend only on $i$, not $x$ or $y$. For, if
$\Gamma_i(x)=W_i(x)-W_{i-1}(x)$ is the set of points at distance $i$ from $x$, then
$$
\begin{array}{lllll}
|\Gamma_1(x)\cap \Gamma_{i-1}(y)| \hspace{-6pt}&=  (q^{e_i}-1)/(q-1),\smallskip
\\
|\Gamma_1(x)\cap \Gamma_i(y)|\hspace{-6pt}&=  (q^{f_i}-1)/(q-1) - (q^{e_i}-1)/(q-1) -1,
\end{array}
$$
and
$$
|\Gamma_1(x)\cap \Gamma_{i + 1}(y)|=  (q^{m}-1)/(q-1) - (q^{f_i}-1)/(q-1) 
$$
(provided also that $i < d)$.  By (g), $|\Gamma_2(x)|=(q^h-q^m)/(q-1)$.  
  By (d) these imply the stated independence.

Counting pairs $ (y,z)$ with $d(x,y)=1=d(y,z)$ and $d(x,z)=2$ yields
$$
|\Gamma_1(x)|\,|\Gamma_2(x)\cap \Gamma_1(y)|=|\Gamma_2(x)|\,|\Gamma_1(x)\cap \Gamma_1(z)|,
$$
whence $(q^m-q)(q^m-q^{f_1})=(q^h-q^m)(q^{e_2}- 1)$.
Equating powers of $q$ yields $1+f_1=m$.  There are then two possibilities:
\begin{itemize}
\item [\rm(i)] $m-1=e_2, ~1=m-f_1=h-m;$  or
\item [\rm(ii)] $m-1=h-m, ~1=m-f_1=e_2.$
\end{itemize}

Suppose (i) holds. Each point is on exactly
$(q^{m-1}-1)/(q-1) = (q^{e_2} - 1)/(q-1) $  ${\mathcal G}$-lines.
Thus, if $d(x,z)=2$,  each of the ${\mathcal G}$-lines on $ z $  contains a point  of the $e_2$-space $W_1(x)\cap W_1(z)$.  Consequently, the graph has diameter $d = 2.$
Moreover, $\Omega $ is a subspace. (For if $x$ and $y$ are distinct points of $\Omega $  but
$\<x,y\>$ is not a ${\mathcal G}$-line, then there is a point $z\in W_1(x)\cap W_1(y)$; then $x$ and $y$ are in the subspace $W_1(z)$, all of whose points are in $\Omega$.)
Now (a) yields $h = n$, so    $m = n- 1  $
and $W_1(x)$  is a hyperplane. Since $y\in W_1(x)$ implies that $x\in W_1(y),$ it follows that  $x\leftrightarrow W_1(x)$ is a symplectic polarity, so (3.1i) holds.

From now on, assume that case (ii) occurs.  Since $e_2 = 1$ there is a unique point
joined to two given points at distance 2. The restriction of the relation ``joined or
equal" to $\Gamma _1(x)$  is thus an equivalence relation, so
 $\Gamma _1(x)$ is a disjoint union of complete graphs, each of size
 $ ( q^{f_1}-q^{e_1})/(q-1) = q(q^{m-2}-1)/(q-1)$.  Since 
 $|\Gamma_1(x)|=  q(q^{m-1}-1)/(q-1) $, this implies that $m-2\mid m-1$, whence $m=3$.
Then $f_1=m-1=2$ (and of course $e_2=1$).

We next determine the sequences $\{e_i\}$,  $\{f_i\}$.  Both are nondecreasing: if
$d(x,y)=i,$ $ d(y,z)=1$ and $d(x,z)=i+1\le d$,  then 
$W_1(x)\cap W_{i-1}(y)\subseteq W_1(x)\cap W_{i }(z)$ and
$W_1(x)\cap W_{i }(y)\subseteq W_1(x)\cap W_{i +1}(z)$.
Also, $e_i<f_i$ since      
$  W_1(x)\cap W_{i-1 }(y) \subset W_1(x)\cap W_{i }(y) $.
If $f_i=3$ for some $i$, then $\Gamma_1(x)\subseteq  W_i(y)$ 
when $d(x,y)=i$,  and so $i=d$;  and conversely
$ f_d = \dim(W_1(x) \cap W_d  (y)) = \dim W_1(x) = 3.$
Thus, $e_i=1$ and $f_i=2$ for $i<d$, while $f_d=3$ and $e_d=1$ or 2.

We will show that ${\mathcal G}$ is a generalized $(2d+1)$-gon or $2d$-gon
(with parameters $q,q$) according as $e_d=1$ or $e_d=2.$
Thus, we must verify axioms (i)--(iii)  given
in Section 2, where $\partial$  was defined. For convenience, we separate the two cases.
 \smallskip

Case $e_d=1$. Since $e_i = 1$ for all $i\ge1$ there is a unique shortest path
joining any two points. Also, a ${\mathcal G}$-line
  $L$ contains a unique point nearest  $x$, unless  $L\subseteq  \Gamma_d(x) $.
  (For, if $y\in  L$ with $d(x,y)=i<d$ minimal, and $u\in W_1(y)\cap W_{i-1}(x)$,
  then $L\ne \<y , u\>  =  W_1(y)\cap W_{i}(x)  $  since $f_i=2$.  
  If $L$ contains a second point in $W_i(x)$ then 
  $L\subseteq U_i(x)\cap \Omega= W_i(x)$ by (f),
whereas  $L\ne W_1(y)\cap W_{i}(x)$.)   
   Thus,  there is a unique shortest path between $x$ and  $L$   if  $\partial (x,L)<2d+1$
   (since then $L \not\subseteq  \Gamma_d(x) $).

Let $L$ and $L' $ be two ${\mathcal G}$-lines.    
If  $L'\not\subseteq  \Gamma_d(x) $ for  all $x\in L$ and
  $L \not \subseteq  \Gamma_d(x') $ for all $x'\in L'$,
  then there is a unique shortest path between
$L$ and $L'$.
(By the preceding paragraph, two shortest paths  would go between points $x_j$ of $L$ and $x'_j$ of $L'$ for $j=1,2$, where $x_1\ne x_2,$ $x_1'\ne x_2'$, and hence produce two shortest paths from $x_1$ to $x'_2$.) 
Suppose $L'\subseteq  \Gamma_d(x) $ for some $x\in L$. Then there is a
unique shortest path from $x$ to each of the $q + 1$ points of $L'$,  no two  such paths
using the same ${\mathcal G}$-line through $x$
(since this would produce a point $y\in W_1(x)$  with $ \partial(y, L') <
2d$ and two shortest paths from $y$ to $L'$).
Then these paths use all   $ q +1$ ${\mathcal G}$-lines
through $x$, and hence $L$ must occur among them.  Thus, $\partial(L,L')=2d$  and a 
unique shortest path again exists  from $L$ to $L'$. Consequently, axioms (i) and  (ii) hold with $n=2d+1$. Since
$f_d = 3$ and $e_d= 1$, so does axiom (iii)   
(using $y$ and any of $q$ $\,{\mathcal G}$-lines on $x$ if $d(x,y)=d$).

 \smallskip 
Case $e_d=2$.   This time, there is a unique shortest path from $x$ to $x'$
unless $x'\in \Gamma_d(x)$. As above, any ${\mathcal G}$-line  $L$
contains a unique point  closest  to $x$, and
there is a unique shortest path from $ x$ to $L$. 
(For, it is not possible for a closest point $y\in L$ to have distance $d$ from $x$, as this would imply that $W_1(y)\cap W_{d-1}(x)$ has dimension  $e_d=2$   and  hence would meet
 $L\subset W_1(y)$ at a point at distance $d-1$  from $x.$)
Finally, let $L $ and $L' $  
be  ${\mathcal G}$-lines with $\partial(L,L')<2d$.   Then only one shortest path can exist between   $L$  and $L'$:   as above, two  would  go  between points $x_j$ of $L $ and $x'_j$ of $L'$ for $j=1,2$, where $x_1\ne x_2,$ $x'_1\ne x'_2$, and hence produce two shortest paths from $x_1$ to $x'_2$.
  Thus, as above axioms (i)-(iii) again hold.

Since $e _2 = 1$, we have $d \ge 3$. The Feit-Higman Theorem (2.3)  now shows
that $d = 3$ and  $e_3 = 2.$

\smallskip
  It remains to prove the parenthetical remark  in (3.1ii).  
 A generalized hexagon with parameters $q,q$  has $|\Omega|=(q^6-1)/(q-1)$ points.
 Since $\Omega =W_3(x)$ is a subspace we have $n=6$.  Since $2=m-1=h-m$ it follows that  
 $W_2(x)$  is a hyperplane and 
 $x\leftrightarrow W_2(x)$ is a symplectic polarity, as required.

\begin{theorem} 

Suppose the hypotheses and conclusions of {\rm (3.1ii)} hold
$($but not necessarily the hypothesis in the parenthetical portion$)$. Then
\begin{itemize}
\item [\rm(i)]If $n = 6,$ then $q$ is even$;$  and
\item [\rm(ii)]otherwise $n = 7$ and $\Omega$  is the set of singular  points of a geometry of
type  $O(7, q).$
\end{itemize}

In either case the embedding of  ${\mathcal G}$  is unique.
\end{theorem} 

We defer the proof to  (A.1iii) in the Appendix.

\section{A reformulation of antiflag transitivity}
Sometimes the following criterion for antiflag transitivity is convenient.

\begin{lemma} 
A subgroup $G$ of $\Gamma L(n,q)$ is antiflag transitive if and only if $G_L^L$ is $2$-transitive for every line $L$.
\end{lemma}

\proof
Suppose $G_x$ has $s$ orbits of hyperplanes on $x$,  $t$ orbits of hyperplanes
not on $x$, and $s' +1$  point-orbits in all. Then $s +t = s' + 1$, and $G_x$ has $s$
orbits of lines through $x$. Each such line-orbit defines at least one point-orbit
other than $\{x\}$. Thus   $t - 1 = s' - s\ge  0$, with equality if and only if $G_{xL}^{L-x}$ 
is transitive for every line $L$ through $x$, as required.
\medskip

From Dickson's list of subgroups of $SL(2, q)$ [5, Chap.~12], it is seen that only
when $q = 4$ is there a 2-transitive subgroup $H$ of $\Gamma L(2, q) $  for which 
$H \cap GL(2, q)$ is not 2-transitive. We deduce the following.

\begin{corollary} 
If $q\ne4$ and  $G\le \Gamma L(n,q)$ is antiflag transitive then so is 
$G\cap G L(n,q),$
and $G_L ^L\cap G L(2,q)\ge SL(2,q)$ for any line $L$.
\end{corollary}
 
\section{The heart of Theorem II}
Suppose  $G\le \Gamma L(n,q)$ is antiflag transitive but not 2-transitive on the points
of $V$. The following lemma incorporates Perin's main idea [20].

\begin{lemma} 
If $x$ is a point$,$  then there is a subspace $W(x)$ $($different from $x$ and $V)$ containing 
$x,$ such that $G_x$ fixes $W(x)$ and is transitive on $V -W(x)$.
\end{lemma}

\proof
A Sylow $p$-subgroup of $G$ fixes a hyperplane $H$ and a point $x\in H$, and is transitive on $V-H$.  Then
$$
W(x)=\bigcap\{H^g\mid g\in G_x\}
$$
is a $G_x$-invariant subspace; $G_x$ is transitive on the pairs $(H^g,y)$ for $g\in G_x, y\notin H^g$, and hence is transitive on $V-W(x)$.  
Finally, $W(x)\ne x$ since $G$ is   not 2-transitive.

\begin{theorem}  
Suppose $G\le \Gamma L(n,q)$ is primitive but not $2$-transitive on points$,$ and is antiflag transitive.  Then $G$ preserves a symplectic polarity$,$   and either
\begin{itemize}
\item[\rm (i)]     $G$ has rank $3$ on points$;$
 or 
 \item[\rm (ii)]   $G$ has rank $4$ on points$,$ $G\le \Gamma Sp(6,q),$
 and $G$ acts on a generalized hexagon with parameters $q,q$   consisting of the points and some of the totally isotropic lines of $V$.
 \end{itemize}
\end{theorem}

The proof involves an iteration of (5.1), followed by (3.1).    Let $d + 1$ denote 
the rank of $G$ in its action on points.

\begin{lemma}  
 There are subspaces 
 $$
 x=W_0(x)\subset W_1(x)\subset W_2(x)\subset  \cdots \subset W_{d-1}(x)\subset  W_d(x)=V 
 $$
 with the properties
\begin{itemize}
\item[\rm (i)]     $G_x$ fixes $W_i(x)$ and is transitive on 
$W_i(x)-W_{i-1}(x)$ for $1\le i\le d;$
 \item[\rm (ii)]    
 if $y \in W_1(x)$ and $0\le i\le d-1,$ then $W_i(y)\subseteq W_{i+1}(x);$
 \item[\rm (iii)]   $W_i(x^g)= W_i(x)^g$ for all $g\in G;$  and
  \item[\rm (iv)]     $d>1$.
 \end{itemize}
\end{lemma}

\proof
Set $W_d(x) = V$ and $W_{d-1}(x) = W(x)$  (cf. (5.1)),  where  $d> 1$  by
 hypothesis.    Since   $W_d(x) - W_{d-1}(x)$  is the largest orbit of $G_x$, certainly 
$W_{d-1}(x^g) = W_{d-1}(x) ^g$     for all $g\in G$.

Now proceed by ``backwards induction". Suppose $W_{j}(x) $ has
been defined for $j = i+1, ... , d,$ and behaves as in (i),  where $i+1<d$;
we need to define $W_i(x)$.
Set  $m_{i+1}= \dim W_{i+1}(x) $.  A Sylow $p$-subgroup $ P$
of $G_x$   fixes a line $L$ on $x$; 
since all 
$P$-orbits on $V -W_{i+1}(x) $  have length at least $q^{m_{i+1} }$,
necessarily  $L\subseteq W_{i+1}(x)$.
If $y\in L-x$ then
 all $P_y$-orbits
on $W_{i+1}(y)-W_{i+1}(x)$ have length at least $q^{m_{i+1}-1}$.
(By primitivity,   $W_{i+1}(y) \ne  W_{i+1}(x)$.)
It follows that  $W_{i+ 1}(x) \cap  W_{i+1}(y)$ is a hyperplane of  $W_{i+ 1}(x) $,
and that $G_{xy}$ is  transitive on $W_{i+1}(y)-W_{i+1}(x)$.  Then 
$$
 W_i(y)=    \bigcap \{W_{i+1}(x)^g\mid  g\in G _y  \} 
$$
is  a subspace of $W_{i+1}(y) $,  and $G_y$ is transitive on $W_{i+1}(y) -W_{i}(y) $.
Then (iii) holds,  since $ G_x$   has only one orbit of size  $|W_{i+1}(x)-W_{i}(x)|$.

This process terminates when $W_0(x) = x$.   Then $W_1(x) -x$ consists of all 
points $y$ for which $\<x, y\>$  is fixed by some Sylow $p$-subgroup of $G.$ 
Now (ii)  follows  from the  definition of $ W_i(y)$. Thus, all parts of (5.3) are proved.

\vspace{4pt}
Let $\mathcal  G$ be the geometry with line set $\{\<x, y\>\mid x\ne y\in  W_1(x)\}$,
 and $\Gamma $  its point graph. By (5.3ii) and induction on $i$, we see that $W_i(x)$ is the set of points at
distance at most $i$ from   $x$  (relative to the metric $d$  in  $\Gamma $). 
  Also,   $G$ is transitive on the pairs $ (x, y)$  with $y \in W_{i+1}(x) - W_i(x)$ for 
each    $i$.
Consequently,
  $\Gamma $  is    metrically regular, and (3.1) applies. 
Since all $W_i(x)$  are subspaces,   (5.2) follows.

\vspace{4pt}

By (3.2), the generalized hexagon in (5.2ii) must be the one associated with
$G_2(q).$    However, as stated in Section 1, we will make the proof of Theorem I,
and most of Theorems II and III, independent of the known existence and uniqueness of
the  $G_2(q)$ hexagon. The required information is easily proved (frequently in the
spirit of other of our arguments), and is collected in the following lemma
  (where $q$ may be even or odd).
   
\begin{lemma}  
 If $G$ is as in {\rm(5.2ii),}  then the following statements hold$:$
\begin{itemize}
\item[\rm (a)]     $G$ has exactly two orbits of t.i.~lines$;$
 \item[\rm (b)]    
$G$ has exactly two  orbits of t.i.~planes$;$
 \item[\rm (c)]   there is a t.i.~plane $E$ such that $G_E^E\ge SL(3, q);$
    
  \item[\rm (d)]     there is an element $t \in G\cap SL(V)$ with $t ^p = 1$ and  $\dim C_V(t) 
  \ge 4;$
  \item[\rm (e)]     $ |G |=   q^6(q^6-1)(q^2-1) c ,$ 
  where $c\,|\, (q-1)e $ if $q=p^e$
  and    $c\,|\, q-1 $ if $G\le GL(V);$
and
    
  \item[\rm (f)]      if $q=2,$ $4$ or $16,$  $r=q+1$  and $R\in \Syl_r( G_x),$  then
$ C_V(R)  $ is a nonsingular $2$-space and $N_G(R)$ is $2$-transitive on $C_V(R)$.
    
 \end{itemize}
\end{lemma} 

\proof
 Since $G_x$ has three point-orbits other than $\{x\}$  (cf.~(4.1)), (a) is clear.
 Clearly, $W_1(x)^G$ is an orbit of $(q^6 - 1)/(q - 1)$ t.i.~planes
(t.i.~using (5.3),  since $W_1(x)^\perp$ is $G_x$-invariant). Let $E$
  be any of the remaining
  $$
  (q^3 + 1)(q^2 + 1)(q + 1) - (q^6 -1)/(q - 1) = q^3(q^3 + 1)
  $$
  t.i.~planes of $V$. If  $L$ is any  ${\mathcal G}$-line, 
  the $q + 1$ t.i.~planes on   $L$
are all of the form $W_1(x)$ for $x \in L$. It follows that $E$ contains no
  ${\mathcal G}$-lines, and for distinct  $y, z\in  E$,   $d(y, z) =   2$;
  let $M = \<y, z\>$ and $x = W_1(y) \cap W_1(z).$ 
Inside $W_1(x) $ there are $q^2$ choices for  $M$, and then there are 
$q$ choices for $E$ on $M$ 
  (any t.i.~plane on $M$
  except $W_1(x)=\<x, y, z\>$). Thus, if    $P \in \Syl_p(G_x)$ then 
 $|P \!: \! P_{ME}|\le q^3$,   so each orbit of $P_{ME}$ on
  the $q^5$ points of  $V - x^\perp$ has length at least $q^2$. Since
 $ E -M$ is fixed by $P_{ME}$, we have $|P \!: \! P_{ME}| = q^3,$ and $P_{ME}$ is transitive on
  $E - M$. This proves (b). Moreover, since $M$ is any line of $E$,  (c) follows from (2.1).
  
  In  (d), let $X< G_E\cap SL(V)$ be a $p$-group inducing  all $(z, \<w, z\>)$-elations (transvections) of   $E$, where
$w\in E - M$. Then $X$ fixes $M$, and hence also the unique point $x$ joined to all of
$M$ by ${\mathcal G}$-lines,
  as well as the unique point $x'$ joined to all of $\<w, z \>$  by ${\mathcal G}$-lines.
   Thus, 
$C_V(X)\supseteq  \<z,x' ,  w , x \>$.

 In  (e), clearly $|G|=(q^3+1)q^3|G_E|$.  
Let $a_1,b_2,a_3,b_1,a_2,b_3$ be the vertices of an ordinary
 hexagon in the point-graph,
with $a_1,a_2,a_3\in E$ and $b _{i+2}=W_1(a_i)\cap W_1(a_{i+1} )$ 
(subscripts mod 3).  Since $W_1(a_1)\cap W_1(b_1)=0$, 
$V=W_1(a_1)\oplus  W_1(b_1)=\<  a_1,b_2,b_3,b_1 , a_2,a_3\>$.
If 
$g\in C_G(E)$  is a $p$-element it follows that $g$ fixes each member of  a basis for $V$.
Then $|G_E|$  divides $|\Gamma L(E)|$  and is divisible by  $|SL(3,q)|$, which implies (e).  
 
In (f),  $R$  cannot fix any point of  $W_1(x)- x$ by (5.3i).  
Then $R$ also cannot fix any point of $x^\perp - x$ 
(if it fixed such a point $y$ it would also fix $W_1(x)\cap W_1(y)$).
Since $R$ fixes a point of $V - x,$  
$ C_V(R)  $ is a nonsingular $2$-space. 
The last part of (f) follows from antiflag transitivity and the Frattini argument.  
    \smallskip 
    
    {\em Remarks.} 1.  Sylow's Theorem and the Frattini argument were standard tools in [19, 12-14], and  will be used  several times below.

\vspace{3pt}

   2. If $q>2$ then $G \cap Sp(6, q)$  is generated by the $G$-conjugates of the group $X$ appearing
in the above proof.

\vspace{3pt}

   3.    If $G \le \Gamma L(n, q)$ is antiflag transitive and primitive on points, then it is
primitive on hyperplanes. For, if $G$ preserves a symplectic polarity, then its 
actions on points and hyperplanes are isomorphic; otherwise, by (5.2), $G$  is
2-transitive on points, and so also on hyperplanes. We will see later (7.1)  that a
stronger result can be obtained by elementary arguments independent of (3.1).
    
\section{The primitive case}
We now begin the inductive part of the proof of Theorems I-III. In order to
avoid identifying $G_2(q)$  during the proof of Theorem I  (cf. Section 1), we
restate the theorems in slightly weaker form.
 
\begin{theorem}  
Let $G \le \Gamma L(n, q), n\ge 2,$ be antiflag transitive. Then one of
the following holds$:$
\begin{itemize}
\item[\rm (i)]    $G \roneq SL(n, q);$
\item[\rm (ii)]    $G$ is $A_7 $ inside $SL(4, 2);$
\item[\rm (iii)]    $G \roneq Sp(n, q);$
\item[\rm (iv)]   
 $G$ is $A_6$  inside $SL(4, 2);$
\item[\rm (v)]   
 $G<\Gamma L(2,4) $
has order $20$  modulo scalars$;$
\item[\rm (vi)]     $G \le \Gamma Sp(6, q) < \Gamma L(6, q),$    and 
$G$ acts as a rank $ 4$ group on
the points of a generalized hexagon  with parameters $q, q,$  whose points and lines
consist of all points and certain  totally isotropic lines for $  Sp(6, q);$
\item[\rm (vii)]    For $q=2$ or $4,$     $G \ron  SL(\frac{1}{2}n, q^2), $  embedded naturally in $\Gamma L(n, q);$
\vspace{2pt}
\item[\rm (viii)]      For $q=2$ or $4,$     $G \ron  Sp(\frac{1}{2}n, q^2), $  embedded naturally in $\Gamma L(n, q);$
  or
 \item[\rm (ix)]   For $q=2$ or $4,$      $G$ is a subgroup of  $ \Gamma Sp(6, q^2),$ itself embedded naturally~in $\Gamma L(12, q),$
such that $G$ acts on a generalized hexagon in $PG(5, q^2)$ as in~{\rm(vi).}
\end{itemize}

 \end{theorem}
  Note that 2-transitive subgroups of     $ \Gamma L(n, q)$
  are automatically antiflag transitive  (Wagner [25, p. 416], or (4.1)).

The theorem will be proved by induction on $n$   in Sections 6, 7. The case
$n = 2$ is omitted, while (2.1) handles $n = 3$.  We therefore assume $n\ge 4.$  By
(4.2), if $q \ne 4$ we may assume that $G \le GL(n, q)=GL(V) $   (compare  (6.1vii-ix)).

In  the remainder of this section we will consider only groups $G$ that are {\em primitive}   on the points of the projective space. Then either (5.2)
applies, or $G$ is 2-transitive. In either case, induction or known results almost
always produce sufficiently large groups of transvections for $G$ to be identified;
case (5.2ii) is exactly (6.1vi), and will be 
considered in the Appendix.

\begin{proposition}
If  {\rm (5.2i)}  holds then either   $G \roneq Sp(n, q) $ or $\,G$ is $A_6$ inside $Sp(4,2)$.%
\end{proposition}
\proof

 We will follow Perin [19]  when possible,  but we include   semilinear groups and the cases
 $Sp(n,2)$  and $Sp(n, 3)$  not dealt with in  [19].
 His method works primarily when $q>4$ and when either $n\ge 6$  or  
  $n=4$ but $q$ is not a Mersenne prime.
  
 If $G$ contains the group of all transvections with a given center, then  $G$ contains all transvections  by transitivity and $G\roneq Sp(n,q)$.  
  
  Assume that  $q>4$, and either $n\ge6$  or $n=4$ and $q$ is not a Mersenne prime.
 We have $| x^\perp -x|=q(q^{n-2}-1)/(q-1)$.
Let $r$ be a primitive  divisor of $q^{n-2}-1$ (see (2.4);
use $r=3$ if  $q^{n-2}-1=8^2-1$ with $n=4,q=8$)  and $R\in \Syl_r(G_x)$. 
Then $r>2$ and   $R< GL(V)$  is  completely reducible,    so    $U=C_V(R)$ is a nonsingular 2-space. 
Moreover,  $N_G(R)^U\ge  SL(2,q)$  (by (4.2), since $G_U=C_G(U)N_G(R)$
by  the Frattini argument),
 while $N_G(R)^{U^\perp}$ is solvable.
Then   $C_G(U^\perp)$ contains $ SL(2,q)= SL(2,q)'$ and hence contains a full transvection group, so  $G\roneq Sp(n,q)$.    Note that the same argument handles the case 
 $G_U^U\ge SL (2,q)  =SL (2,4)$.

It remains  to consider the possibility that either $q\le4$ or  that  $n=4$ 
and $q$ is   a Mersenne prime. 

Let $ x$ and $y$ be distinct points of the t.i.~line $L$. There
is a Sylow $p$-subgroup $P$ of $G$  fixing $x$ and $L$, and transitive on $V - x^\perp$. Then all
orbits of $P_y$ on $V - x^\perp$  have length at least $q^{n-1}/q$, 
so $P_y$ is transitive on $y^\perp -x^\perp$.
Since $G_y$  is already transitive on $y^\perp/y$ by (5.2i), it is  antiflag transitive there.

By our inductive hypothesis concerning (6.1), $K =G_y^{y^\perp/y}$   satisfies one of the
following conditions:%

\begin{itemize}
\item[\rm ($\alpha$)]    $K \roneq Sp(n-2, q);$
\item[\rm ($\beta$)]    $K=A_6 ,n-2=4,q=2;$
\item[\rm ($\gamma$)]    $K$ acts on a generalized hexagon as in (6.1vi), $n-2=6;$
\item[\rm ($\delta$)]    For $q=2$ or 4,   $Sp(\frac{1}{2}(n -2),q^2)\nor K
\le \Gamma  Sp(\frac{1}{2}(n -2),q^2)$  with $\frac{1}{2}(n -2)$ even$;$
\item[\rm ($\epsilon$)]      For $q=2$ or 4,      $K < \Gamma Sp(6, q^2)$ acts on a generalized hexagon over $GF(q^2)$ as in (6.1vi),
$n-2=12;$
or
\item[\rm ($\zeta$)]     $K<\Gamma L(2,4)$ has order $20$  modulo scalars, $n-2=2$.
\end{itemize}

In particular, if $q$ is odd then $K\roneq Sp(n-2,q)$.

  If $Q=O_p( Sp(V)_y )$ and $T$ is the group of transvections in $Q$,  then $Q/T$ and $y^\perp/y$ are naturally $\Gamma Sp(V)_y$-isomorphic projective
modules   (via $uT\to [V,uT]/y$ for $u\in Q$).
Moreover,  $T$ is the Frattini subgroup of $Q$
if $q$ is odd, while $Q$ is naturally an $O(n-1,q)$-space if $q$ is even.

The case $G\le Sp(n,q) = Sp(4,2) \cong S_6$ is easily handled and so will be excluded.   Note that
$|G\cap Q|=4$ 
if $G=Sp(4,2)' \cong A_6$.

If  $G\cap Q\not \le T$ we will show that $G\ge T$ and hence $G$ contains $Sp(n,q)$.
 Let $r$ be a primitive divisor of $ q^{n-2}-1$ and $R\in \Syl_r (G_y\cap Sp(n,q))$  (using 
 $r=3$ if  $q^{n-2}-1=8^2-1$ when $n=4,q=8,$ or
 $r=7$ when $n=8,q=2$).   
 The $R$-invariant subgroup $W=[G\cap Q,R]$  projects onto a subspace $WT/T$ of 
  the $GF(q)$-space $Q/T $; 
  in view of the action of $R$ on $y^\perp/y$ and hence on $Q/T$, 
  we have     $WT/T=Q/T $.
If  $q$ is odd it follows that $W$ contains the Frattini subgroup $T$ of $Q$, so that $G>T$.
If  $q$ is even then $W$ is a  nonsingular hyperplane of the orthogonal $GF(q)$-space $Q$.  
 If $G$ does not contain $T$ then each element of $G_y$ leaves  the hyperplane $W$ invariant, while acting antiflag transitively on $ y^\perp/y$ and hence on the 1-spaces of the  orthogonal space 
 $W$,
 so  we are in  case ($\z$).
Then  
     $G_{yL}^L$  is $\Z_4$  by (4.1).  Since $W$ is elementary abelian, $|W^L|\le2$
 and $|C_W(L)|\ge8$.  If $L'$ is a second t.i.~line containing $y$ then
 we obtain the contradiction
$4\le |C_W(L)\cap C_W(L')|=| C_W(\<L,L'\>) | = | C_W(y^\perp ) |= 1$.%

   If    $ n=4$   then $q^3\le |G_y|_p=|G\cap Q||K|_ p\le 
|G\cap Q|qe _p$  where $q=p^e$, so that $|G\cap Q | >|T|$ and we have seen 
   that  $G\roneq Sp(4,q)$.
This takes care of dimension $n=4$, including  ($\z$).
From now on $n\ge6$ and $q\le 4$.

If $n=6$  then the same argument yields 
$q^5\le |G\cap Q| |Sp(4,q)|_pe_p = |G\cap Q| q^4e_p$, so $G\cap Q\ne1$.  
We have already handled the cases   $G\cap Q\not \le T$ and $G\cap Q \ge T$.    
 It remains to eliminate the possibility   $1\ne G\cap Q< T$, where $p|e$ and hence $q=4$.
 Since the above inequality shows that $(\a)$ holds,  
if $E\subset y^\perp$ is a nonsingular 2-space  then some
$g\in G_{yE} \cap Sp(n,4)$  induces an element of order $3$ on 
$E/y$ and hence acts in that manner on a nonsingular 
2-space $D\subset \<y,E\>$, fixing a point $z \in D$.
Some $h\in G$ satisfies $z^h=y\in D^h$, and then $(g^h)^{D^h}$
acts nontrivially on $(G\cap T)^{D^h}$,
which contradicts the assumption  $1<| G\cap Q |<4$.

Now $n>6$.
If $q=3$ then $K\roneq Sp(n-2,3)$.  Let $r$ be a primitive   divisor of $3^{n-4}-1$ and $R\in \Syl_r(G_y)$.
Then     $U=C_V(R)$ is a   nonsingular 4-space.  
Since  $R$ is a Sylow subgroup of the stabilizer of two  
  perpendicular points of $U$ and of the stabilizer of two
  non-perpendicular points of $U$,  by the Frattini argument
   $N_G(R)^U$  has rank 3  and hence
 contains $Sp(4,3)$  by induction.
Also  $N_G(R)^{U^\perp}$ is solvable
(lying in  $\Gamma L(1,3 ^{n-4})$).
Then  $C_G(U^\perp)$ contains  transvection groups and 
 $G\roneq Sp(n,3)$.

Now  $q=2$ or 4.
 In   $(\a)$ let  $r$  be a primitive divisor of $q^{(n-2)-2} - 1$
(use $r = 7$ if $q^{(n-2)-2} - 1= 2^{(10-2)-2} - 1$),
 in  $(\g)$ let $r=q+1$,
in $(\d)$ let  $r$  be a primitive divisor of $ (q^2) ^{\half(n-2)-2}-1$, 
and in   $(\e)$   let $r=q^2+1$. 
 Let  $R \in  \Syl_r(G_{yL})$. 
Then   $C_{y^\perp/y}(R)$ is a nonsingular 2-space over $GF(q)$  in $(\a)$ 
and  $(\g)$  (cf.   $(5.4f)$)
or over $GF(q^2)$  in
  $(\d)$ and $(\e)$, so $ U = C_V(R)$  is 
nonsingular of dimension 4 or 6. As above,  by the Frattini argument
   $N_G(R)^U$  has rank 3 and hence 
    contains   $  A_6 $, $Sp(4,q)$ or $Sp(6,q)$ by induction.
    Also    $N_G(R)^{U^\perp}$  is solvable or is a subgroup of
    $\Gamma Sp(2,q^2)$ in $(\c)$ or   of
    $\Gamma Sp(2,q^4)$   in $(\e)$.  Then 
    $N_G(R)$  has a subgroup  $N$  inducing the identity on $U^\perp$  and
 $A_6$, $Sp(4,q)$ or $Sp(6,q)$  on $U.$  
 In the last two cases we obtain   $G \roneq Sp(n,q)$ as usual; in the $A_6$ case 
  a $G$-conjugate of $N$ meets $Q$ in a subgroup of size 4 and hence
    $G\cap Q\not\le T$, which was handled above. 
    
\vspace{2pt}
  The next primitive case is  Theorem I. 

\begin{proposition}
If  $G\le \Gamma L(n,q)$  $(n\ge3)$  
is $2$-transitive on points$,$   then either $G\ge SL(n,q)$ or $G$ is $A_7$ inside $SL(4,2)$.
\end{proposition} 
 
\proof
In view of Wagner [25, Theorem 4], we may assume that $n \ge6$. We recall the
following additional facts from Wagner [25, pp. 414, 416]: $G$ is 2-transitive on hyperplanes,
 and
if  $H$ is a hyperplane, then $G_H^H$ is antiflag transitive.

Once again,  we will run through the possibilities provided by induction for
$G_H^H$  and, dually,  $G_x^{V/x}$.   If either is 2-transitive, then
by induction  $G$  is  
transitive on  complete flags (i.e., maximal increasing
sequences of subspaces, one of each dimension), and the
result follows from Wagner [25, Theorem 3] or Higman [8, Theorem~1]; so suppose not.  Let $H$ be a hyperplane and $x\in H$.

Suppose $G_x^{V/x}$  is primitive and hence  is contained in $\Gamma Sp(n-1 ,q)$ by (5.2).  
Then $G_{xH}$ fixes a line $ \Delta $ on $ H$ and $x$.
By (4.1),  $(G_{H }^{H })^\Delta_{ \Delta}$ is 2-transitive, so $G_{xH}<  G_{\Delta H}$ and  $G_H^H$ is imprimitive. 

Thus,  we may assume that $K=G_x^{V/x}$ is imprimitive.
By Theorem  III,  $q=2$  or~4,
 $n-1\ge 6-1$ is even and $K\le \Gamma L(\frac{1}{2}(n-1), q^2)$  behaves   
 as follows:

\begin{itemize}
 
\item[\rm ($\alpha$)]    $K \ron  SL(\frac{1}{2}(n-1), q^2);$
\vspace{2pt}
\item[\rm ($\beta$)]    $K \ron  Sp (\frac{1}{2}(n-1),q^2);$
or
\vspace{2pt}
\item[\rm ($\gamma$)]      $K < \Gamma Sp(6, q^2)$ acts on a generalized hexagon over $GF(q^2)$ as in (6.1vi),
$n-1=12.$
\end{itemize}
 
Let $r$ be a primitive divisor of $(q^2)^{\half (n-1)-1  }-1$ in  ($\alpha$)
or   of  $(q^2)^{\half (n-1) -2 }-1$ in  ($\beta$), and  let  
  $r=q^2+1$ in   ($\gamma$).
(Use $r=7$ if $(q^2)^{\half (n-1)-1  }-1 = (2^2)^3-1$ in  ($\alpha$).)
 Let  $R\in \Syl_r(G_{x})$.
Then   $\dim C_V(R)$ is $1+2$ in ($\alpha$)  
and $1+4$ in ($\beta$)  and  ($\gamma$)
(using (5.4f) in   ($\gamma$)).
By the Frattini argument, $N_G(R)$ is 2-transitive on $U=C_V(R)$,  inducing at least $SL(U)$  by induction. 
Moreover, $N_G(R)^{[V,R]} $ is solvable, except perhaps in ($\gamma$)
with $N_G(R)^{[V,R]}\le \Gamma L(2,q^{4})$.
As usual, $C_G([V,R])^U \ge SL(U)$, so  $G$ contains 
a full transvection group and $G\ge SL(V)$, which contradicts the behavior of $K$.

\smallskip

\em {Now  \rm(5.2), (3.2), (6.2) \em   and \rm(6.3)  \em complete the inductive step in 
\rm(6.1)\,\em   when $G$ is primitive on points.}\rm
 
 \smallskip
 Having dealt with the primitive case, we record an elementary corollary for
use in the next section.
\begin{lemma}    
Suppose $G$ is as in 
{\rm(6.1)} and is primitive on points.~If $F\le G$~with $F$ antiflag transitive and 
$|G \!: \!  F|$ a power of $p , $ then $F$ is also primitive on points.
\end{lemma}

\proof 
 Let $P\in {\Syl}_p(G_x)$. 
Then $P$ fixes a  unique line $L$  on $x$.  
(In case (6.1vi),  by (5.3i)  the $p$-parts of the nontrivial orbit lengths of $G_x$ are $q$, $q^3$ and $q^5$.)

Clearly $G=P\hspace{.5pt}F$ and $P\cap F\in {\Syl}_p(F_x)$. 
If $F$ is imprimitive then, by (6.1vii-ix),
  there is a unique line containing $x$ fixed by $F$,  and it is also the unique line fixed by  $P\cap F$;
  this line   must be $L$. 
Thus,   $G_x=P\hspace{.5pt}F_x$ fixes $L$, contradicting (4.1)  and the  primitivity of $G$.
  
\section{The imprimitive case; completion of the proof}

Continuing our proof of (6.1), we now turn to the case of an antiflag transitive
subgroup $G$ of $\Gamma L(n, q)$ that is {\em imprimitive on points.}  The method here is
entirely different from that of Sections 5, 6; we build a new projective space on
which $G$ continues to act antiflag transitively.

If  $\Delta$  is a nontrivial imprimitivity block for the action of $G$ on points, then $\Delta$ 
is the set of points of a subspace. (For, 
every hyperplane  of $\<\Delta\>$ does not contain some point of $\Delta$.  Then
$G _\Delta $  is transitive on the hyperplanes of
$\<\Delta\>$, hence on its points, and thus $\Delta$  must contain all points of 
$\<\Delta\>$.) We usually
identify $\Delta $ with  $\<\Delta\>$.   Set $\delta = \dim  \Delta$ and
 $W\cap \Delta^G=\{\Delta'\in \Delta^G \mid \Delta'\subseteq W\}$
  for any subspace $W$.

By Remark 3 at the end of Section 5, $G$ is
also imprimitive on hyperplanes, and a block of imprimitivity consists of all
hyperplanes containing a subspace $\Sigma$. The next result (independent of the aforementioned
Remark) shows that there is a close connection between blocks of
points and hyperplanes. It is due to Orchel [16], and simplifies and improves a
result in an earlier version of this paper.

\begin{lemma}[Orchel]  
Let $\Delta $  be a block of imprimitivity for $G$ acting on points$,$
and $\delta = dim  \Delta .$ Let $H$ be a  hyperplane$,$  and  let $\Sigma$ be the union of the members~of  $\Delta^G$ contained in $H$. Then $\Sigma$ is a subspace of dimension $n - \delta$  partitioned by $\Delta^G\cap\Sigma ,$ and the set of hyperplanes containing $\Sigma$ is a block of imprimitivity for $G$ acting on hyperplanes.

\end{lemma}
\proof
We have $|\Delta^G|=(q^n-1)/(q^\delta-1)$.   Then $(q^{n-1}-1)/(q-1)= 
|H\cap \Delta^G| \cdot\break(q^{\d}-1)/(q-1)
+  \big((q^{n}-1)/(q^\d-1) -  {|H\cap \Delta^G| \big)   (q^{\d -1 }-1)/(q-1)}
$, so    $|H\cap \Delta^G|=(q^{n-\delta}-1)/(q^\delta-1)$.  
The union $\Sigma$ of the members of $H\cap \Delta^G$ has cardinality 
$ ( q^{n-\delta}-1)/(q-1)$.  
 
Let $P\in \Syl_p( G_H)$. Then $P$ is transitive on $V- H$, and hence on 
$  \Delta^G-(H\cap \Delta^G)$.
Let $\Sigma' $ be a subspace of $H$ of dimension $n-\delta$ fixed by $P$.  If $\Sigma'\cap \Delta'\ne0$ for one (and hence all) 
$\Delta'\in \Delta^G-(H\cap \Delta^G)$,
then $|\Sigma'|\ge| \Delta^G-(H\cap \Delta^G)|= q^{n-\delta}$, which is false; so $\Sigma \supseteq \Sigma' $, and comparing cardinalities shows that $\Sigma= \Sigma '$   is a subspace.
Moreover, if $\Sigma\cap \Delta'\ne0$ for $\Delta' \in \Delta^G$, then 
$\Delta' \subset H$ and hence $\Delta' \subseteq \Sigma$,   so
  $\Sigma$ is partitioned by $\Delta^G\cap\Sigma .$

Now, if $H'$ is any hyperplane containing $\Sigma$, then   
$H'$ contains all $\Delta^g\subset H$, so $H'\cap \Delta^G$ contains
$H\cap \Delta^G$
and hence has union containing $\Sigma$.
As any element of $G$ sending $H$ to $H '$ sends  
 $H \cap \Delta^G$ to  $H'\cap \Delta^G$,
$G_\Sigma$ is transitive on the set of such hyperplanes $H'$.  This proves the lemma.
\smallskip\smallskip

{\em Notation.} 
 Let $\Delta$ be a minimal proper block of imprimitivity, and define $\Sigma$ as in (7.1).  
Let ${\mathcal L}$ be the set of all intersections of members of $\Sigma^G$.

\begin{lemma} 
If $n>2\delta$ then ${\mathcal L}$
is the lattice of subspaces of a projective space $PG(n/\delta -1,q^\delta)$ on which $G$ acts as an antiflag transitive collineation group.
\end{lemma}
\proof
By (7.1), if $W\in  {\mathcal L}$ then
$W\cap \Delta^G$ partitions $W$.  If $W=\<\Delta_1,\dots,\Delta_k\>$ 
with $\Delta_i\in \Delta^G$ and $k$ minimal, then $\dim W=k\delta$ and 
$|W\cap \Delta^G|=(q^{k\delta}-1)/(q^\delta-1)$.
Call $W$ a Point, Line, or Plane if $k=1,2$ or 3, respectively.  Then two Points are on a unique 
Line (containing $q^\delta+1$ Points), 
and three Points not on a Line are  in a 
unique Plane (containing $q^{2\delta}+q^\delta+1$ Points).
The  Veblen and  Young axioms  [24] imply that  ${\mathcal L}$  is a projective space.

By (7.1),  $H\cap \Delta^G=\Sigma\cap \Delta^G$,
 and $G_H$ is transitive on the $q^{n-\delta}$ Points not in $\Sigma$.  
Thus, $G$ acts antiflag transitively on ${\mathcal L}$.

\smallskip\smallskip
\textsc {Definition}.  Let ${\mathfrak  A}$ denote the set of all cosets of members of 
${\mathcal L}$.  (Since $0\in  {\mathcal L}$, all vectors of $V$ are in ${\mathfrak  A}$.)
 
\begin{lemma} 
If $n>2\delta$ then ${\mathfrak  A}$
 is the lattice of subspaces of $AG(n/\delta,q^\delta)$.
\end{lemma}
\proof
 Form  ${\mathfrak A}\dot\cup {\mathcal L}$ by attaching ${\mathcal L} $   ``at infinity" as follows: adjoin $U\in {\mathcal L}$ to $W+v$ if $U\subseteq W\in {\mathcal L}$.
 Thus,  ${\mathfrak  A}\dot\cup {\mathcal L}$ will have two types of ``points" (vectors 
and members of $\Delta^G$), and two types of  ``lines" (cosets of members of $\Delta^G$, and  Lines of  ${\mathcal L}$).
If $\<\Delta, \Delta'\>$ is a Line of ${\mathcal L} $, then it and any vector determine a translation plane of order $q^\delta$ in a standard manner [4, p.~133]; $\<\Delta, \Delta'\>$ 
plays the role of line at infinity. By
(7.2),  ${\mathfrak  A}\dot\cup {\mathcal L}$ satisfies the Veblen and Young axioms, and hence is $PG(n/\delta,q^\delta)$.
This proves the lemma.

\begin{lemma} 
If $n=2\delta$ then ${\mathfrak  A}$
 is  $AG(2,q^\delta)$.
\end{lemma}
\proof
As above,  ${\mathfrak  A}$ is an affine translation plane of order $q^\delta$. But here $\Delta^G$ is merely its
line at infinity, so proving that    ${\mathfrak  A}$   is desarguesian will be more difficult.
 We will use standard  properties of  collineations of finite projective planes   [4, Chap.~4].  Using dimensions, $V=\Delta\oplus \Delta'$ for distinct 
 $\Delta, \Delta'\in \Delta^G$.

Let $x \in \Delta$ and $P \in \Syl_p(G_x)$. The group $E = C_P(\Delta)$   
consists of all elations of  ${\mathfrak  A}$ 
with axis $\Delta$; it is semiregular on 
the set $\Delta^G - \{\Delta\}$ of lines $\ne \Delta$  of  ${\mathfrak  A}$  through the
point 0 of  ${\mathfrak  A}$, and  ${\mathfrak  A}$  is desarguesian if $|E |= q^\d,$ by (2.2). We may thus assume that  $|E| < q^\d$   and aim at a contradiction.

Let $H \supset \Delta $ be a hyperplane fixed by $P$, so
$H\cap  \Delta^G =\{\Delta\}$ by (7.1).
Since $P$ is transitive on $V-H$ and hence on $\Delta^G - \{\Delta\}$,
$P_{\Delta'}$ is transitive on $\Delta'-H$ (since $\Delta'$ is a block).
Now $G_\Delta$ is transitive on the pairs $(x,\Delta')$ with $x\in \Delta$ and $\Delta'\in \Delta^G-\{\Delta\}$, so $G_{\Delta\Delta'}^\Delta$ is transitive
and hence $G_{\Delta\Delta'}^{\Delta '}$  is  antiflag transitive 
since $P^{\Delta '-H}$ is transitive.
 Moreover, $G_\Delta=P\cdot G_{\Delta\Delta'}$ since $P$ is transitive on 
 $\Delta^G-\{\Delta\}$.  Then  $G_{\Delta }^\Delta=P^\Delta G_{\Delta\Delta'}^\Delta$;
 since $G_{\Delta}^\Delta$  is primitive by the minimality of $\Delta$, 
 $G_{\Delta\Delta'}^\Delta$
   is primitive by (6.4).

 We claim that $C_G(\Delta)_{\Delta'}=1$.  For,
 $C_G(\Delta)\noreq G_\Delta,$  where $G_\Delta$ is transitive on $\Delta^G - \{\Delta\}$ and $C_G(\Delta)_{\Delta'}$ consists of homologies of  ${\mathfrak  A}$ with axis  $\Delta$.
 Thus, if $C_G(\Delta)_{\Delta'}\ne1$, then this holds for every 
 $\Delta'\in \Delta^G-\{\Delta\}$.  Then in the action of $C_G(\Delta)$ on $\Delta^G-\{\Delta\}$,
 the stabilizer of any two points is trivial, but the stabilizer of   any point is nontrivial.  This implies that $C_G(\Delta)$ acts as a transitive Frobenius group on $\Delta^G-\{\Delta\}$, with kernel   $E$ of order $q^\delta$, contrary to assumption.
 
 It follows that $C_G(\Delta)=E$, and 
 $|G_\Delta^{\,\Delta} \!: \! G_{\Delta\Delta'}^{\,\Delta} |=
 (|G_\Delta|/|E|)/|G_{\Delta\Delta'}|=  q^\delta/|E|$  is a power of $p$.
 
 Suppose $q$ is odd. By (4.2) we may assume that  $G\le GL(n,q)$.  By induction,
 both $G_{\Delta\Delta'}^{\,\Delta}$ and $G_{\Delta }^{\,\Delta}$ have normal subgroups $SL(\delta,q)$ or $Sp(\delta,q)$ or a group as in (5.4).    The known orders (cf.~(5.4e))
 do not allow for distinct subgroups of one of these types to have index a power of $p$ in one another (since  $G\le GL(n,q)$).  
 It follows that 
$G_{\Delta\Delta'}^\Delta=G_{\Delta }^\Delta$, and $q^\delta/|E| =1$,  contrary to assumption.

Consequently, $q$ is even.  Since $G_{\Delta\Delta'}$ has even order it has an involution $t$.  
 Then $t$ is a Baer involution (since   it fixes $\Delta $  and $\Delta'$), and   
 $\dim C_\Delta(t)=\half \delta$
and $|C_E(t)|\le   q^{  \delta/2}$ (since $C_E(t)$ acts on  the Baer subplane for $t$).
Induction for $G_{\Delta\Delta'}^{\,\Delta}$, together with this restriction on 
 involutions in $G_{\Delta\Delta'}$ (cf. (5.4d)), imply that either
($\alpha$) $\delta=2$, or  ($\beta$) $\delta=4, $ $q=2$, 
$G_{\Delta\Delta'}^{\,\Delta}=A_6$ or $A_7$.
\smallskip 

($\alpha$)   The argument used for $q$ odd applies, unless $q = 4$, $G_\Delta^{\,\Delta}= SL(2,4).2$
and  $ G_{\Delta\Delta'}^{\,\Delta} = SL(2,4)$
(modulo scalars).  Here, $4^2/|E|=q^\delta/|E|=
|G_\Delta^{\,\Delta} \!: \! G_{\Delta\Delta'}^{\,\Delta} |=2$,~so    
$G_{\Delta\Delta' }\cong SL(2,4) $   centralizes $E$. 
 Choosing $t$ in this $SL(2,4)$ contradicts $|C_E(t) | \le 4$.

\smallskip 

($\beta$) In this case,  $|\Delta^G|=2^4+1$  and $|G_\Delta^{\,\Delta} \!: \!  G_{\Delta\Delta'}^{\,\Delta}|=8$  or 2, corresponding to (6.1ii,iv), so $|E|=2 $ or 8.  
  Since $G_{\Delta\Delta'}^{\,\Delta}\ge A_6$,
  \vspace{2pt}
   the argument in $(\alpha)$ yields
   $|E|=2$.  Then $G_{\Delta\Delta'}$ 
fixes both members of    $\Delta'{}^E $, so 
  $G_{\Delta\Delta'}$ fixes $k\ge 3$ points and we obtain a Steiner system
$S(2,k,17) $, which is impossible. 

\smallskip 
This completes the proof of (7.4).
\smallskip \smallskip

\emph{Proof of} (6.1).
Each translation $v\to v+c$ permutes the members of 
${\mathfrak  A}$, sending each hyperplane  of ${\mathfrak  A}$ to 
itself  or a disjoint hyperplane.  Thus,  these form  the group of translations of the affine space  ${\mathfrak  A}$.   
The corresponding group of scalar transformations acts homogeneously on $V$, and hence is uniquely determined 
up to $GL(V)$-conjugacy.
Then the group $G^+$ of all
collineations of ${\mathfrak  A}$   induced by elements of $\Gamma L(n, q)$  is   $\Gamma L(n/\delta , q^\d)$.

In particular, $(G_\Delta^+) ^\Delta  = \Gamma L(1, q^\d) $ has  order   
$(q^\d-1)\d e$, where $q=p^e$.  Since this group is antiflag transitive, 
$q^{\d-1}$ divides $\d e$, whence $q=\d=2$
or  
$e=\delta=2, $  $ q=4$. 

 Thus, $G$ is an antiflag transitive subgroup of $\Gamma L(\half n,q^2)$,
where $\Gamma L(\half n,q^2)$ is embedded naturally in $\Gamma L(V)$.
Moreover, $G$ acts primitively on the set $\Delta^G$ of points of  
$PG(n/\delta -1,q^\delta)$. For otherwise, there
is an imprimitivity block $\Lambda \supset \Delta$, and $G$ lies in 
$\Gamma L(\half (\half n),(q^2)^2)$.   Then $q=2$, $|\Lambda|=16$ and $G_\Lambda^\Lambda$ lies in 
$\Gamma L(1, 16)$, which we have just seen is not antiflag transitive.

This primitivity and (6.1i,iii,vi) produce (6.1vii-ix),    finishing our proof of (6.1).

\smallskip
 \smallskip
  {\em Remark.}
  {\em Examples  of}   (6.1vii-ix) {\em  occur.}
 For, let $F=GF(q^2)\supset K=GF(q)$ with $q=2$ or 4,  and
 $V=V(\half n ,q^2)=Fv\oplus W$ with $n\ge4$ even and $W$ an  $F$-hyperplane.
  Let $b \in F-K$
and $\sigma\in \Aut( F)$  of~order $\log_2q^2=q$.
Then the $K$-hyperplanes containing $W$ are
$Kv\oplus W$ and $Kb^{\sigma^i}v\oplus W$, $0\le i<q$, where
  $\<\sigma\>$  fixes the first of these and is transitive on the remaining ones.
 Since $SL(\half n,q^2 )$ contains a subgroup of order $q+1$ 
 transitive on the 1-dimensional $K$-subspaces of $Fv$, 
  $SL(\half n,q^2 )\<\sigma \> $ is antiflag transitive.
 The symplectic and $G_2$ cases are similar.
  Moreover, {\em any antiflag transitive instance of} (6.1vii-ix)
 \emph{contains one of the   groups
generated by $SL(\half n,q^2),$ $Sp(\half n,q^2)$ or  $G_2(q^2)$
together with a group of $q$ field automorphisms.}

\smallskip
\smallskip 
Now the proofs of  (6.1) and Theorem I are complete.  Moreover, for Theorems II and
III, we  only have to identify the groups occurring in (6.1vi) -- the hexagon  ${\mathcal G} $  is already known  to be both unique and correctly embedded, by (5.2) and (3.2).~It is known that the group of automorphisms 
of  ${\mathcal G}$ induced by elements of $Sp(6, q)$ is $G_2(q)$; this is stated 
in Tits [22, (11.3)] and  proved in Tits [23,  (5.9)]. We observe, independent of this,
that $G \cap Sp(6, q) = G_2(q)$:   in view of $G_2(q) \le \Aut( {\mathcal G} )$  and  (A.6iii),
 if $S $  is 
the group of scalar transformations of $V$  then  $|GS\cap  GL(6, q) |= |G_2(q)S|$ and $G_2(q)\cap S = 1$.

\section{Corollaries}
In this section we give some consequences of Theorems I-III.

The {\em  affine group}  $A \Gamma L(n, q)$  is defined as the group
$$
\{v\to v^g+c\mid g\in \Gamma L(n,q) , \, c\in V\}
 = T\semi \Gamma L(n,q) $$
 of all collineations of the affine space $AG(n, q)$ based on $V$, an $n$-space over 
 $GF(q)$. ($T$ denotes the translation group.)

\begin{proposition} 
Let $G \le  A\Gamma L(n, q), n\ge3,$ be transitive on ordered non-collinear triples of points of 
$AG(n, q).$ Then $G = T \semi G_0,$  where $T$ is the translation 
  group$,$ and $G_0\roneq  SL(n, q)$ or $G_0$ is $A_7$ $($with $n = 4,$  $ q = 2).$
\end{proposition}
\proof
 The hypothesis implies that $G_0$  (the stabilizer of 0) is projectively one of the
groups of Theorem I; it remains only to show that $G$ contains $T$. If not, then
$G \cap T = 1$ (since $G_0$  is transitive on points), 
and so $|G |\le |\Gamma L(n, q)|$
since $|A\Gamma L(n, q)  |=| \Gamma L(n, q) |  |T|$.
But then $|G \!: \!  G_0| = q^n $ contradicts   $|\Gamma L(n, q) \!: \! G_0|\le  (q - 1)e$  (resp. 
 $|\Gamma L (n, q) \!: \!  G_0| = 8$) if $G_0\ge  SL(n, q),$ $ q = p^e$  (resp. $G_0 = A_7).$

\begin{corollary} 
The only proper $3$-transitive subgroup  of $A\Gamma L(n,2)$ is
$V_{16}\semi A_7$ when $n=4$.
\end{corollary}

This corollary improves various results in the literature (for example Cameron [3,
Theorem 1]); and also Jordan's theorem (Wielandt [26, (9.9)]):

\begin{corollary} 
A normal subgroup $N$of a $3$-transitive group $G$ is $2$-transitive$,$
unless it is elementary abelian of order $2^n$ and either $G=N\semi GL(n, 2)$
or $n=4$ and $G=N\semi A_7$.
\end{corollary}

From results of Perin [19] and Kantor [12], we deduce the following

\begin{proposition}
Suppose $G \le \Gamma L(n, q )$ is transitive on the $j$-subspaces of $PG(n-1,q)$ for some $j$ with $2\le j\le n-2.$
Then $G$ is transitive on the $i$-subspaces for all $i$ with $1\le i\le n-1,$ and one of the following  occurs$:$
\begin{itemize}
\item[\rm (i)]   $G\roneq SL(n,q);$
\item[\rm (ii)]  $G$ is $A_7$ inside $GL(4,2);$  or
\item[\rm (iii)]  $G$ is $\Gamma L(1,2^5)$ inside $GL(5,2)$.
\end{itemize}
\end{proposition}

{\em Remark.}   A ``$t$-$(v, k, \lambda)$  design in a finite vector space" is a collection of
$k$-subspaces or  ``blocks" in a $v$-space, any $t$-space being contained in precisely
$\lambda$  blocks.   No nontrivial examples are known with $t\ge 2$; and (8.4) shows that
none can be constructed by the analogue of the familiar construction of $t$-designs
from  $t$-homogeneous groups (Dembowski [4, (2.4.4)]).
\smallskip

To motivate the next result, we sketch the deduction of Perin's Theorem [20]
(mentioned in Section 1) from Theorem II. Suppose $G \le \Gamma L(n, q), n\ge 4,$ and
suppose $G$ acts as a primitive rank 3 group of even order on the points of
$PG(n - 1, q)$. For a point $x$,   $G_x$ has three orbits on points, and hence three
orbits on hyperplanes. If $G$ is antiflag transitive, then $G \le\Gamma Sp(n, q)$ by
Theorem II (and indeed $G$ is known). Otherwise, $G_x$  is transitive on the hyperplanes
through $x$, and so also on the lines through $x$, in contradiction to Kantor [12].

\begin{proposition} 
Suppose $G \le \Gamma L(n, q),$ $ n\ge 4,$  and $G$ acts as a primitive
rank $4$ group on the points of $PG(n - 1, q).$   Then either $q = 2,$ $ 3, $ $4 ,$ 
$5$ or $ 9,$ or
$G \roneq G_2(q),$  $q$ even$, $ embedded naturally in  $ \Gamma Sp(6, q).$

\end{proposition} 

\proof 
 By Theorem II, we may assume that $G$ is not antiflag transitive;
by the previous argument and Kantor [12], we may assume it is not transitive
on incident point-hyperplane pairs. Thus, of the four $G_x$-orbits on  hyperplanes,
two consist of hyperplanes containing  $x$. Then $G_x$  has two orbits on lines containing
$x$.  There are thus two $G$-orbits on lines, with $G_x$ transitive on the lines of
each orbit which pass through $x$. Consequently, $G_L^{\,L}$  is transitive for each line $L.$
 
 Since  $G_x$  has three orbits on points different from $x$, it follows that, for suitable
$L$ and  $M$  from different line-orbits,  $G_{xL}^{L-x}$   is transitive while 
$G_{xM}^{M-x }$ has
two orbits. Thus, $G_L^L$  is 2-transitive while $G_{M}^{M }$ has rank 3.
  But,  using Dickson's list of subgroups of $PSL(2, q)$   [5, Chap.~12], we see that $P\Gamma L(2, q)$  has a rank 3 subgroup only if $q = 2, 3, 4,5$ or 9.

\begin{proposition} 
Let $G$ be an irreducible subgroup of   $P \Gamma L(n, q),$ $ n \ge 4.$
Suppose $G_x$  is  transitive on the lines through $x,$ for some point  $x$.  Then $G$ is $2$-transitive on points $($and Theorem I applies$).$
\end{proposition} 

\proof  By Kantor [12], it is enough to show that $G$ is transitive on points.
So let $X =  x^G$ and assume $X$ is not the set of all points. If  $L$ is a line and
$L \cap X \ne  0,$  then $\ell = | L \cap X|$  is independent of $ L$, and $1 < \ell < q - 1$. If
$\dim W = m$ and $W\cap X\ne 0,$ then  $|W\cap X | = 1 + (\ell- 1)(q^{m-1}  - 1)/(q - 1).$

There is an $(n - 2)$-space $U$ disjoint from $X$ (for otherwise the
hyperplane sections of $X$ would be the blocks of a  design  having the same $b, r, \lambda$ as $PG(n -1 ,q)$ and hence  $|X|=v=b$). The
hyperplanes containing $U$  partition $X$ into sets of cardinality    
$k=1+ (\ell - 1)
(q^{n-2} - 1)/(q - 1)  $; so $k$ divides $ |X|= 1+ (\ell- 1)(q^{n-1}-1)/(q-1) $
and hence also $ (\ell-1)q^{n-2}$. 
Then $(q-1)-(l-1)\equiv 0$ (mod $k$). 
Since $k > (q^{n-2} - 1) /(q - 1) > q$, we have  $\ell = q.$  
But then the complement of $X$ contains one or all points
of each line,  and so is a hyperplane fixed by $G$, contradicting irreducibility.

 \bigskip
\centerline{II.  THEOREMS IV AND V}

\section{The geometry of primitive antiflag transitive groups}

The proof of Theorem IV occupies Sections 9-11. The present section contains
notation  and the analogue of (5.3). The primitive case is concluded in Section 10;
there the method is different from that of Section 6.  Unlike Theorems I-III, the
primitive case   does not depend on the imprimitive one. Finally,  Section 11 corresponds
to Section 7.

The symplectic case is covered by Theorems II and III; so {\em we will exclude the
case  $ G \le \Gamma Sp(2 m , q)  $ for   the remainder of the proof.} Also, in view of the isomorphism between the $Sp(2 m, q)$ and $O(2   m+1, q)$ geometries when $q$ is even,
{\em we will  also exclude the case $G \le \Gamma O(2 m  +1, q), $  $q$ even}.   Thus, the geometry is
associated with a nondegenerate sesquilinear form.

In the proof,  $\Omega$  denotes the set of  t.i.~or t.s.~points of the appropriate classical geometry, defined on a
vector space $V$ over $GF(q)$. (This assumption involves  
{\em a slight change of notation
in the unitary case}:  
$G$  will be a subgroup of $\Gamma U(n, q^{ {1/2}})$.
 This may lead to the impression of minor discrepancies between the statement of
Theorem IV and parts of Sections 9-11:   the notation for the name of the group will remain the same as in Section 2, only the meaning of   ``$q$'' will change.) 

In general our convention is to refer only to 1-spaces in $\Omega$, though there will be situations where other 1-spaces will be mentioned.    Thus, 
in general we identify a subspace with the set of members of $\Omega$ it contains;  some care is needed when dealing with anisotropic subspaces.   Similarly, in general
if $S$
is a subset of  $\Omega$, then $S^\perp$ is the set of points of  $\Omega$ collinear with (i.e.,~perpendicular to) every point of $S.$ 
The subspace $0$ plays the role of $\emptyset$, so $0^\perp =  \Omega$.
 This convention has   odd-looking consequences, such as: 
a t.i.~or t.s.~subspace $W$ is maximal if and only if $W^\perp = W.$ 
(However, if  $W$ is nonsingular and if no point is collinear with every point of $W$,
then $W^\perp$ will denote an anisotropic vector subspace.)
The notation $\<X\>$ usually refers to  a vector subspace, not just a set of points;
the meaning will be clear from the context.
The dimension of a t.i.~or t.s.~subspace is its vector space dimension
(cf.~Section 2), and the {\em rank}  $ r$ of the geometry is the maximal such dimension.

 \vspace{2pt}
We begin with two preliminary lemmas.

\begin{lemma} 
 \hspace{-5pt}
There do not exist subspaces $T,$ $W\!$ with 
$T\cup T^\perp =W^\perp \!$~and~$T,$$\,\,T^\perp\ne W^\perp$.
\end{lemma}
\proof
If $T\cup T^\perp =W^\perp $ then $T\cap T^\perp =(W^\perp)^\perp =W$.
Let $t_1\in T-W$ and $t_2\in T^\perp-W$, and observe that a point of $\<t_1,t_2\>-\{t_1,t_2\}$  is not in  $T\cup T^\perp$.

\begin{lemma}   

Suppose $T,W$ are t.i.~or  t.s.~subspaces with $\dim T=i-1,$
$\dim W=i,$ and $T\subset W$.  Then $|T^\perp - W^\perp |=q^{2r-i    +   c},$
where $c\ge -1$ depends on the type of $V$ but not on $r=\rank(V)$ or $i,$ and is given in the following table.
$$
 \begin{tabular}{ccccccc}
\hspace{-6pt}
\!\!\!Type\,of  $V\,$	&\hspace{-6pt}$\!Sp(2r,q)$   &\hspace{-6pt}$\!O^+(2r,q)\!$ &\hspace{-6pt}
$\!\!O(2r\!+\!1,q)$ &\hspace{-6pt}$\!O^-(2r\!+\!2,q)$ 
&\hspace{-6pt}$\!U(2r,q^{ {1/2}})$   &\hspace{-6pt}$\!U(2r\!+\!1,q^{ {1/2}})\!\!$%
 \raisebox{2.4ex}{\hspace{-1pt}}\raisebox{-1.1ex}{\hspace{-1pt}}
 \\
 \hline
$ c$&$ 0$         &$ -1 $&$0 $ &$  1 $&$ -\half $ &$\half$  
\raisebox{2.7ex}{\hspace{-1pt}}\raisebox{-1.4ex}{\hspace{-1pt}}
\
\raisebox{2.7ex}{\hspace{-1pt}}\raisebox{-1.4ex}{\hspace{-1pt}}
\end{tabular}
$$
\end{lemma}
\proof
For $i = 1,$  $|T^\perp - W^\perp|=|\Omega - W^\perp|$ is the number of  points not 
perpendicular to the point $W$, and is easily computed.
For $i\ge2$, $T^\perp/T$  has rank  $r-i+1$ and the same type as $V$;
each of its points outside $W^\perp/T$ corresponds to a coset
 (containing $q^{i-1}$ points) of $T$ outside $W^\perp$.
 
 \smallskip\smallskip
 Throughout the rest of this section and the   next, $G$ will be assumed to act
{\em antiflag transitively} on the geometry and 
{\em  primitively on the set $\Omega$ of points}.
 Let $d + 1$ denote the rank of $G$   on points.

\begin{lemma}  
   For each point $x$ there is a chain of $G_x$-invariant  subspaces 
 $
 0=W_{-1}(x)\subset
 {x=W_0(x)}
 \subset W_1(x)\subset   \cdots \subset    W_d(x)=V 
 $
 with the following  properties$:$
\begin{itemize}
\item[\rm (i)]    $W_i(x)^\perp =W_{d-i-1}(x)$  $($in particular$, $
$W_i(x)$ is t.i.\!  or t.s.~if and only if $i\le \frac{1}{2}(d-1));$
 \item[\rm (ii)]      $G_x$   is transitive on  $W_i(x)-W_{i-1}(x)$ for each $i;$
 \item[\rm (iii)]    
  if $y \in W_1(x)$ and $0\le i\le d-1,$ then $W_i(y)\subseteq W_{i+1}(x);$

  \item[\rm (iv)]      $W_i(x^g)= W_i(x)^g$ for all $i,$ $x,$ $g;$ and
  \item[\rm (v)]   $W_1(x)\cap W_{1}(y)$ is a hyperplane of $W_1(x)$ if $y\in  W_1(x) -x $ 
  and $d\ge 4$.
 \end{itemize}
\end{lemma}

\proof  Let $L$ be a line on $x$ fixed by some $P \in \Syl_p (G_x).$
For $y\in L-x$, all $P_y$-orbits on $V-x^\perp$ have length at least $q^{(2r-1+ c)-1}$
by (9.2), so $P_y$ is transitive on 
$y^\perp - L^\perp$ (again by (9.2)).  Set $W_1(y)=\< L^g \mid g\in G_y\>$. Then 
$y\in L\subset y^\perp$, so $y\in W_1(y)\subset y^\perp$.  Moreover,  
$$
W_1(y)^\perp=\bigcap \{ (L^\perp)^g \mid g\in G_y\}\!,
$$
and $G_y$ is transitive on  $y^\perp- W_1(y)^\perp$. 
(In particular,  $W_1(y)^\perp$   does not depend on $x$: it is the unique 
$G_y$-invariant subspace $U$ of $y^\perp$ such that $G_y$ is transitive on 
$y^\perp- U$.)
  Define  $W_1(y^g)=W_1(y)^g$ for all $g\in G$.

If $W_1(x)^\perp = x$, we are finished (and $d=2)$.  So suppose $W_1(x)^\perp \ne  x$.
Then  $W_1(x)\cup W_1(x)^\perp \ne  x^\perp$, by (9.1).  Since $G_x$ is 
 transitive on  $x^\perp- W_1(x)^\perp$, it follows that $W_1(x)\subseteq W_1(x)^\perp$, that is, $W_1(x)$ is t.i.~or t.s.  
(in the characteristic 2 orthogonal case $W_1(x)$ is t.s.~since it is totally isotropic and spanned by t.s.~subspaces).
Also, $G_x$ is transitive on $W_1(x)-x$.  (For, $W_1(x)$ is naturally isomorphic to the dual space of $V/W_1(x)^\perp$.   Now $G_x$ has two orbits on 
the points  of  $V/W_1(x)^\perp$, namely those in $x ^\perp/W_1(x)^\perp$  and 
those not in $x ^\perp/W_1(x)^\perp$; 
so it has two orbits on the points of $W_1(x)$,
namely $x$ and $W_1(x)-x$.) 

Now proceed by induction, assuming that  $1\le i\le \frac{1}{2}(d-1)$ and that subspaces  $W_j(x) $ and  $W_{d-j-1}(x)$ 
have been defined for $-1\le j\le i$ satisfying (i)-(iv).
Set   $m= \dim W_i(x)$.  
By (ii) and (9.2), the $P$-orbits on $V-W_i(x)^\perp$ have length at least $q^{2r-m+ c}$,
and hence the $P_y$-orbits have length at least $q^{(2r-m+ c)-1}$.
We may assume that $m\ne r$, since otherwise we are finished.
Again by (9.2),  $q^{2r-m+c-1}\ge  q^{m}>|W_i(x) |$.
  As above,
$W_i(y) \subseteq W_i(x)^\perp$ and   $\< W_i(y) , W_i(x)\>$ is t.i.~or t.s., where 
$ W_i(x) \ne W_i(y)$  by primitivity.   Since $ P_y$ acts on 
$  W_i(y)^\perp - \< W_i(y) , W_i(x)\>^\perp$ with   orbit  lengths at least 
$q^{2r-(m+1)+c}$, (9.2) implies that 
$W_i(y)$ is a hyperplane of  $\< W_i(y) , W_i(x)\>$ and $P_y$ is transitive on 
$W_i(y)^\perp -\< W_i(y) , W_i(x)\>^\perp $.

Set $W_{i + 1}(y) = \<W_i(y),W_i(x)^g \mid g \in  G_y\> \subseteq W_i(y) ^\perp$. Then $G_y$ fixes $W_{i+1}(y)$
 and is transitive on  both 
  $W_i(y)^\perp -W_{i+1}(y)^\perp $ and    $W_{i+1}(y) -W_{i}(y)$. (For, as before,   $W_{i+1}(y)$ is naturally  isomorphic to the dual space of 
 $V/W_{i+1}(y)^\perp$, and  $G_y$ has   exactly $i+1$ point-orbits  
 $(W_{j-1}(y)^\perp /W_{i+1}(y)^\perp    )-  (W_{j}(y)^\perp/W_{i+1}(y)^\perp )$,   $0\le j\le i$,
 on  $V/W_{i+1}(y)^\perp$, while  acting on $i+1$ subsets  
 $W_j(y) -W_{j-1}(y) $  of
 $W_{i+1}(y)$.)
 
 Now (i,ii,iv) hold, while  (iii) follows from the definition of $W_{i+1}(y)$ if $i\le \half (d-1)$,
 and from $W_{d-i-1}(y)\supseteq W_{d-(i+1) -1 }(x)$ if 
 $i > \half (d-1)$. This completes the inductive~step.  
 
 Finally, (v) was proved in our argument when $i = 1$, since $m\ne r$ in that  case.

 \smallskip\smallskip
 
\textsc {Definition}.   The geometry ${\mathcal G}$  consists 
of the points of $\Omega$,
  together with those lines 
  ({\em${\mathcal G}$-lines}) joining $x$ to points of
   $W_1(x)$ for all $x \in   \Omega$.  The point graph of  ${\mathcal G}$ 
  is~$\Gamma.$
  By (4.1), if $y$ is a point of  a ${\mathcal G}$-line $L$ then $L\subseteq W_1(y)$.

\begin{lemma}   { \rm (i)}    $\Gamma$  is metrically regular.
\begin{itemize}
 \item[\rm (ii)]     $d\le4.$
 \item[\rm (iii)]     If $V$ has type $O(2 r +1,q),$   
 then the conclusions of {Theorem IV} hold.
  \item[\rm (iv)]     If $d = 2$ then the conclusions of {Theorem IV}  hold.
 \end{itemize}
\end{lemma}

\proof 
(i) This follows from (9.3ii-iv).

(ii) If  $d\ge 5$ then $W_2(x)$ is t.i.~or t.s., 
 and hence  satisfies axiom  (g) in Section 3, so (3.1) yields a contradiction.

(iii)  Recall  that $W_2(x)$ is either $x^\perp$ or t.s., and hence
$ |W_2(x) |= (q^h - 1)/(q - 1)  $ for some $h$. If $d = 2$ then $G$ has rank 3 on points,
 and Kantor-Liebler
[14, (1.3)] applies, since $q$ is odd. If  $d = 3$ then (3.1) and (3.2) show that ${\mathcal G}$ 
 is the generalized hexagon associated with $G_2(q)$, embedded naturally in $V$ of type $O(7,q)$.  Then $G\roneq G_2(q)$ as at the end of  Section 7.

(iv)   Use
 Kantor-Liebler [14, (1.3),\,(6.1)]  (since we have excluded the symplectic case).

\smallskip \smallskip

{\em Notation.}  $ e_2 = e$ and  $f_1 = f $ are defined as in Section 3; $W(x) =W_1(x) $,  and
$m = \dim W(x).$

\begin{lemma}  
 $d = 4$  is impossible.
 \end{lemma}  
 \proof
 If $ d = 4$  then the chain in (9.3) is
 $$
 0\subset x\subset W(x)\subset W(x)^\perp \subset x^\perp \subset V,
 $$
the differences being orbits of $G_x$.  By (9.3v), $f=m-1$.  Let $N_{r-m}$
denote the number of   points of $W(x)^\perp/W(x)$.
Then $|W(x)^\perp -W(x) |=  q^m N_{r-m}$, as in the  proof of (9.2).  As in Section 3, 
a count of pairs $(y,z)$ with $d(x,y)=d(y,z)=1$, $d(x,z) = 2$,  yields
$$
(q^m-q)(q^m-q^{m-1})= q^m N_{r-m}(q-1)(q^e -1).
$$
  Thus, $e$ divides $m-1$.
 
 Since  $W(x) \ne  W(y)$  for  $x \ne y,$  $ W(x)$ is not a clique. Let $y, z \in W(x) $  be nonadjacent points. Then 
  \begin{align*}
 e &=  \dim W(y) \cap W(z)
\\
&\ge \dim W(x) \cap W(y) \cap W(z)
\\
 &\ge   m-2,
  \end{align*}
 since $W(x) \cap W(y)  $  and $W(x) \cap   W(z)$ are hyperplanes of $W(x)$.
 Now 
 $q^{m-1}-1=N_{r-m}(q^e-1)>q^e-1$, $e\mid m-1$
 and $e\ge m-2$ force $m\le 3$.  Clearly $m>2$, since $\Gamma$ is connected.  Thus, $m=3$  and $N_{r-m}=(q^2-1)/(q^e-1)$.  Then $e=1$ and 
 $|W_2(x)|=|W_1(x) |  + q^mN_{r-m}= (q^5-1)/(q-1)$, in contradiction to (3.1).
 
\section{The case $d=3$}
In this section we continue the proof of Theorem IV in the primitive case.
By Section 9 we may assume that $d = 3 $ and $V$ is not of type  $O(2r + 1, q).$
The chain of subspaces in (9.3) is now
$$
0 \subset x \subset W(x) \subset  x^\perp \subset V,
$$
with $W(x)$ maximal t.i.~or t.s.  Set  $k=|x^\perp -x|$ 
and $v_i=(q^i-1)/(q-1)$.

\begin{lemma}  
$V$ has type $O^+(2r, q)  $  with $r=  4, 5$ or  $6, $  while $f =r- 2$ and
$e =2.$
 \end{lemma}  
 \proof
 As usual, count the pairs $(y, z)$  with $d(x, y) = d(y, z) = 1, $  $d(x, z) =2, $  this time obtaining
$$
(v_r-1)(v_r-v_f) = (k-(v_r-1))v_e.
$$
In particular, $k\le (v_r-1)(v_r-v_f + 1) \le v_r(v_r-1)$.
However,  $k$ is easily computed for each type, and the types  $O^-(2r+2,q)$ and $U(2r+1,q^{1/2})$
fail to satisfy this inequality. 
Moreover, in the case $U(2r,q^{1/2})$, we have $k=  q  (q^{r-1}-1)(q^{r-3/2}+1)/(q-1)$, whence 
$$
v_r-v_f=q^{r-3/2} v_e ,
$$
and   $f = r- 3/2$, which is absurd.

Thus, $V$ has type $O^+(2r,q)$.  This time,
$$
k-v_r+1 = q^{r-1}(q^{r-1}-1)/(q-1)=q^{r-2}(v_r -1),
$$
so
$$
q^{r-2} v_e =v_r-v_f,
$$
whence $f=r-2$, $e=r-f=2$.  
Since $W(x)\cap W(y)\supseteq \<x , y\> $ for $y\in  W(x)-x$,  $f\ge2$, so $r\ge4$.

Let $y, z$ be nonadjacent vertices in  $W(x)$. Then
  \begin{align*}
\hspace{106pt}
 2 &=  e \ge   \dim  W(x) \cap W(y) \cap W(z)
 \\
&\ge 2(r-2)-r \hspace{170pt}(*)
\\
 & = r-4,
  \end{align*}
whence $r\le 6$, as required.

\begin{lemma}  
$r = 4.$
 \end{lemma}  
 \proof
 Suppose $r= 5$ or 6. 
  If $\<x,y\>$ is a  ${\mathcal G}$-line then $\dim W(x)\cap W(y)=f=r-2>2$, so there is a
 point $z\in W(x)\cap W(y)-\<x,y\>$.  
 Call the span of three noncollinear but pairwise
adjacent points a {\em special plane};
note that all lines of a special plane   belong to
 ${\mathcal G}$ (since $\<x,y\>\subset W(z)$). 
  Then 
 $$
 |W(x)\cap W(y)\ -\<x,y\>| = (q^f-q^2)/(q-1),
 $$
 so $\<x,y\>$   lies in exactly $(q^{f-2}-1)/(q-1)$ special planes.   
 If $f=r-2=3$, this number is 1, so the number of special planes is 
 $$
 v_5(q^4+1)\cdot v_4\cdot 1/(q^2+q+1)(q+1),
 $$
 which is not an integer.  So $r=6$   and $f=4$. 
 
 In this case, we will show that the ${\mathcal G}$-lines and special planes that pass through
$x$ form a generalized pentagon with parameters $q, q$, contradicting the Feit-Higman Theorem (2.3).
 
 Any special plane through $x$  contains $q + 1  $ \ ${\mathcal G}$-lines through $x$, and any such
${\mathcal G}$-line lies in $(q^{f-2}-1)/(q-1)=q+1$ special planes.  If $\<x,y\>$ and  
$\<x,z\>$ are  ${\mathcal G}$-lines
through $x$ not contained in a special plane, tightness in the inequalities  $(*)$
shows that $W(x)\cap W(y)\cap W(z)$  is a  ${\mathcal G}$-line through $x$, 
the unique such  ${\mathcal G}$-line
lying in special planes with both $\<x,y\>$ and   $\<x,z\>$.    Now elementary counting  
verifies axioms (i)-(iii) for a generalized pentagon in Section 2,  which yields the desired contradiction.

\smallskip\smallskip

{\em There are several ways to handle the case} $r = 4$. One is to show that ${\mathcal G}$ is a ``dual
polar space'' (of type $O(7, q))$ in the sense of Cameron [3]; another is to quote
transitivity results in Kantor-Liebler [14, Sect. 5]. The method used here involves
triality, a concept which we now briefly discuss;
we will see that triality is involved in the embedding appearing  in  Theorem IV(iii).
 We refer to [22] for further discussion of triality.

Let  $ {\mathcal P}$  be the set of points of the geometry of type $O^+(8, q),$   ${\mathcal L}$ the setá of lines,
and   ${\mathcal M}_1$  and $ {\mathcal M}_2$ the two families of 
{\em solids}  (maximal t.s.~subspaces); thus, any t.s.~plane
lies in a unique member of each family. More generally, two solids lie in the
same family if and only if their intersection has even dimension. The geometry
admits a ``triality automorphism'' $\tau$ mapping  $ {\mathcal L} \to {\mathcal L}$ and 
${\mathcal P}\to {\mathcal M}_1\to {\mathcal M}_2\to {\mathcal P}$ and preserving the natural incidence between 
$ {\mathcal P}  \cup {\mathcal M}_1\cup {\mathcal M}_2$ and $ {\mathcal L}$ (defined by inclusion or reverse inclusion).  Also, $\tau$ preserves the ``incidence'' on 
$ {\mathcal P}  \cup {\mathcal M}_1\cup {\mathcal M}_2$, in which a solid is incident with a point contained in it, and two solids are incident if they meet in a plane.  This  ``automorphism''
induces an automorphism of $ P\Omega^+(8,q)$.

\smallskip 

Before continuing with the proof, we outline the way in which the examples of
Theorem IV(iii) arise. Let  $v$ be a nonsingular  vector, so that  
$v^\perp \cap {\mathcal P} $ carries a
geometry of type  $O(7, q).$ If $M_i \in {\mathcal M}_i\, (i = 1, 2),$   then
$v^\perp \cap M_i $ is a plane, contained
in a unique member $M_i ^*$ of ${\mathcal M}_{3-i};  $   thus $v$ induces bijections between
${\mathcal M}_1, $  ${\mathcal M}_2$ and the set of 
planes (maximal t.s.~subspaces) of $v^\perp \cap {\mathcal P} $.  
These
bijections are invariant under $G =\Omega^+(8,q)_v$, which acts transitively on each set.  Now apply triality: $G^\tau$ 
is an irreducible subgroup of  $\Omega^+(8,q)$, transitive on  ${\mathcal M}_2^\tau={\mathcal P} $,   
and preserving a ``geometry" on  ${\mathcal P} $  isomorphic to the dual polar space of t.s.~planes of  $v^\perp \cap{\mathcal P} $. 
(Strictly, here and below, in place of $G^\tau$ we use the inverse image in 
$\Omega^+(8,q)$ of $(G/Z)^\tau$, where $Z=Z(\Omega^+(8,q))$.)
$G$ is transitive on disjoint pairs of t.s.~planes of $v^\perp \cap{\mathcal P} $, and hence on disjoint pairs of elements of ${\mathcal M}_2$; hence $G^\tau$ is transitive on nonperpendicular members of ${\mathcal P} $, that is, antiflag transitive.
Note that $G^\tau$ and $G^{\tau^{-1}}$  lie in different conjugacy classes in 
$\Omega^+(8,q)$.  Note also that $ G^\tau = \Omega(7,q)$  only if $q$ is even; for  $q$ odd, 
$ G^\tau $  contains the element $-1\in \Omega^+(8,q)$.

The process can be continued one further time. If $w\in V$ is a nonsingular vector, then $ (G^\tau)_w$ acts transitively (and even antiflag transitively) on $w^\perp \cap  {\mathcal P}$, preserving
a geometry that is the $G_2(q)$ hexagon, naturally embedded.

\smallskip \smallskip

We return to the proof. There are
$(q^4-1)(q^3 + 1)/(q-1)=(q+1)(q^2+1)(q^3+1)$ points, and equally many subspaces $W(x)$.  
Since $f=e=2$, $\dim W(x)\cap W(y)=2$ or 0 for $x\ne y$, and so
 all subspaces $W(x)$ belong to the same family; without loss of generality, $\{W(x)\mid x\in {\mathcal P} \} =
{\mathcal M}_1$.   

Call $M\in {\mathcal M}_2$ {\em special}  if it contains a ${\mathcal G}$-line
$L$.    If   $x$ is a point of $M-L$,
then $\<x,L\>$ is contained in  a unique member  $W(y)$ of ${\mathcal M}_1$,
 and $M\cap W(y)= \<x,L\>$.  Since ${\mathcal G}$ has no triangles
 (as $f=2$), we have $y\in L$, and   $\<x,y\>$ is a ${\mathcal G}$-line.
Thus, the points and ${\mathcal G}$-lines in $M$ form a  generalized quadrangle
using all points of $M$.
Then   $x\leftrightarrow M\cap W(x)$ is a symplectic polarity of $M$
whose absolute lines are the  ${\mathcal G}$-lines in $M$;
so the quadrangle is of type $Sp(4,q)$.

Let $\Lambda$ be the set of special solids, so   $\Lambda^\tau$ is a  set of points.  We claim that, {\em if $U$ is a solid$,$ then  $U\cap  \Lambda^\tau$   is a 
t.s.~plane.}
For, 
 $U^{\tau ^{-1}}\cap  \Lambda$ is the set of special solids incident with 
 $U^{\tau ^{-1}}$, which is either a point     $x$ or  a solid $W(x)$,
 so  $U^{\tau ^{-1}}\cap  \Lambda$
    is  the set of special solids $M$ such that 
  $M\cap W(x)$ is a plane of  $W(x)$  containing $x$.
  That set is  the set of special solids incident with the incident pair $x, W(x)$,
  and hence
 the set of  solids incident with the  pair $x, W(x)$.  
  Its $\tau$-image is the set of points incident with an incident pair $x^\tau,W(x)^\tau$ of 
  solids, and hence is  a t.s.~plane  of $V$.
\smallskip

The following result now identifies $\Lambda  ^{\tau}$ (and hence  $\Lambda$).
\begin{theorem}
Let $\Phi$ be a subset of $\Omega,$ the point set of a geometry of type $O^+(2r,q),$ $ r\ge3$.
Suppose that$,$ for every t.s.~$r$-space $U$ of $\Omega,$  ~
$\Phi\cap U$ is an $(r-1)$-space. 
 Then   $\Phi=\Omega\cap  v^\perp$ for some nonsingular vector~$v$.
 
\end{theorem}

\proof
We treat first the case $r = 3$. Identify $\Omega$  (the Klein quadric) with the
set of   lines of   $PG(3, q)$. Then a t.s.~plane of $\Omega$  is either the set of lines on a point
or  the set of lines in a plane; and a line of $\Omega$  is the set of lines in  a plane $E$ and
on a point   $x\in  E$. Thus, under this identification, $\Phi$  is a  set of lines of $PG(3, q)$
having the property that the members of $\Phi$  on a point $x$  all lie in a plane $E$, while
those in a plane $E$ all contain a point $x$.   Then $x \leftrightarrow E$ is a symplectic polarity,
and  $ \Phi$ is its set of t.i.~lines.  A symplectic polarity of $PG(3, q)$  can be
identified with a  point  $v$ outside the Klein quadric $\Omega$, its t.i.~lines corresponding
to points of $   \Omega  \cap   v^\perp  $.

For  $r > 3$ we use induction. If  $x, y$  are  nonperpendicular points  of  $\Phi$,  then $\Omega\cap \<x,y\>^\perp =\Omega'$ is of type $O^+(2r-2,q)$.
We claim that $\Phi\cap \<x,y\>^\perp =\Phi'$ satisfies the conditions of the theorem in $\Omega'\,$(with $r-1$ replacing $r$).
If $U$ is a t.s.~$(r-1)$-space in $\Omega'$, then $\<x,U\>$  is a t.s.~$r$-space, and 
$ \Phi\cap \<x,U\>$ is an 
$(r-1)$-space containing $x$ by hypothesis; so $\Phi  \cap U=\Phi ' \cap U$  is an $(r-2)$-space.
By induction, $\Phi\cap \<x,y\>^\perp=\Omega \cap  \<x,y,v\>^\perp$ for a nonsingular vector $v\in   \<x,y\>^\perp$.

For $a,b\in \Phi$   distinct and perpendicular,     
 a t.s.~$r$-space  $U$   containing $a$  and  $b$ produces     
  an $(r-1)$-space $\Phi\cap U$, so $\<a,b\>\subseteq \Phi\cap U \subseteq \Phi$. 

In particular, $\Phi \cap x^\perp  \supseteq \< \<x,b \>\mid b\in \Phi \cap \<x,y \>^\perp  \>  = \Omega \cap \<x,v\>^\perp$
and $\Phi \cap y^\perp  \supseteq   \Omega \cap \<y,v\>^\perp$.  Every point of 
$\Omega \cap v ^\perp$ lies on a line meeting $\Omega \cap \<x,v\>^\perp$ and
$\Omega \cap \<y,v\>^\perp$ in different points, so $\Omega \cap v^\perp\subseteq \Phi$.
Finally, $ \Phi \subseteq \Omega \cap v^\perp $, since each point of 
$\Omega-v^\perp $ is in a  t.s. $r$-space  properly containing 
a t.s.   $(r-1)$-space of $v^\perp$, and hence cannot lie in $\Phi$.

\smallskip

{\em Remark.}  The theorem fails for $r = 2,$ $q>3$:  $\Omega$ is a ruled quadric 
(a   $(q+1)\times (q+1)$ square lattice), and there are $(q+1)!$  sets $\Phi$ satisfying the hypothesis of (10.3), only $(q+1)q(q-1)$ of which are conics.

\smallskip \smallskip
{\em Completion of the proof of the primitive case of Theorem IV.}
It remains to identify $G$.  Let $H$ be the group induced by $G^\tau$ on the $O(7,q)$ 
geometry $\Lambda^\tau \subset v^\perp$.  Then $H$ is transitive and has  rank 4 
  on the set of t.s.~planes contained in $\Lambda^\tau$.
  (For,  $H$ has rank 4 on ${\mathcal P} ^\tau= {\mathcal M} _1$ and hence 
  on the set of planes $W(x)\cap \Lambda^\tau$.)
We will use the action on these planes to show that
$H$ contains $\Omega(v^\perp)$.

If $E$ is a plane, then $H_E$ is transitive on the t.s.~planes meeting $E$ in a line, so
$H_E^E$ is line-transitive.  Also, $H_E$ is transitive on the 
$q^6$  t.s.~planes disjoint from $E$.
Since any point outside $E$ lies on $q^3$ such planes, every point-orbit outside $E$
of a Sylow $p$-subgroup $P$ of $H_E$ has length divisible by $q^3$.  Let $L$ be a line of $E$ fixed by $P$.
Since $L$ only lies in $q$ t.s.~planes $E'\ne E$, each of which has $q^2$ points outside $E$,  it follows that $P_{E'}$ is transitive on $E'-E$.  
By (2.1),  $H_E^E\ge SL(E)$.
If $x$ is any point of $E$   then  $C_H(x)_E$ is transitive on $E/x$.
Since $E$ can be any t.s.~plane of $v^\perp$ on $x$, it follows that $C_H(x)$ is transitive on $ x^\perp/x$.

Let $Q$ denote the centralizer of both $x$ and $x^\perp/x$ in $\Omega(7,q)$.
We have $q^6\le |H|_p=|H\cap Q||H_x^{x^\perp /x} |_p \le |H\cap Q|q^4e$ since 
$H_x^{x^\perp /x}\le \Gamma O(5,q)$, so
  $H\cap Q \ne 1$.
But $Q$ is elementary abelian of order $q^5$, and is $C_H(x)$-isomorphic to
$x^\perp /x$.
Then $C_H(x)$ acts irreducibly on $Q$, and hence $H\cap Q =Q$.
If $h\in H$ and $x^h\notin x^\perp$, then $H\ge \<Q,Q^h\> =\Omega(7,q)$.

This completes the primitive case of Theorem IV.

\section{The imprimitive case }

Throughout this section (which corresponds roughly to Section 7),  $G$ satisfies
the hypotheses of Theorem IV and is imprimitive on points.  We are assuming
that $V$ has rank $r\ge3$.  

Let $\Delta$ be a proper block of imprimitivity for $G$ on $\Omega$.  Then   $G_\Delta^{\,\Delta}$ is transitive,
while $G_x=G_{x\Delta}$ is transitive on $V-x^\perp$ for $x\in \Delta$.  Thus, $\Delta\subseteq \Delta^\perp$, and $G_\Delta$ is transitive on $V-\Delta^\perp$.
Then $\<\Delta\>$ is t.i.~or t.s., and (by the duality between $V/\<\Delta\>^\perp$
and $\<\Delta\>$) \,$G_\Delta$ is transitive on $\<\Delta\>$.  Thus, 
{\em$ \Delta = \<\Delta\>$
is a t.i.~or t.s.~subspace   and $G_\Delta^\Delta$ is antiflag transitive.}

From now on, $\Delta$ will be a minimal proper block of imprimitivity and $x\in \Delta$.
Set $\d=\dim \Delta$.    

\begin{lemma} 
$  x^\perp \cap  \Delta^G= \Delta^\perp \cap  \Delta^G$  partitions  $ \Delta^\perp$.
\end{lemma}
\proof    
This is clear if $\Delta $ is a maximal t.i. or t.s. subspace, so assume that 
$\d<r$.
Suppose   
  $\Delta'\in \Delta^G$,  $\Delta'\cap \Delta^\perp\ne0$  and
   $\Delta' \not\subseteq \Delta^\perp$.   Let $y\in \Delta- \Delta'^\perp$. 
Since $G_y$ is transitive 
on $\Omega-y^\perp $ it is transitive on $\Delta^G-(  y^\perp\cap \Delta^G  ) $, so every member of $\Delta^G-( y^\perp\cap \Delta^G ) $ meets $\Delta^\perp-\Delta$.  
Since $ | \Delta'-y^\perp| = q^{\delta-1}$,  
$$
|\Delta^\perp-\Delta|\ge |\Delta^G-(  y^\perp\cap \Delta^G )|
=|\Omega-y^\perp|/q^{\d-1}.
$$
A check of each classical geometry (computing $ |\Delta^\perp-\Delta|$
 as in   (9.2)) shows that this inequality does not hold except in the case
 $O^+(2 r,2) $.  
 
 Consider that case.  We  will use a different inequality that is stronger in that case.  
Since $G_\Delta$ is transitive on $\Delta$,   every member of 
 $\Delta^G-(   \Delta^\perp  \cap \Delta^G  ) $ arises for some $y$ as above,
 so $G_\Delta$ is transitive on this set.
 If $k=\dim ( \Delta^\perp\cap \Delta' )$ then  
\begin{align*}
 |\Delta^\perp -\Delta| &\ge    |  \Delta^G-(   \Delta^\perp  \cap \Delta^G  )  |  
  (2^k-1) 
   =   \{ |\Omega -\Delta^\perp |  /  \! ( 2^\d-2^k )\}(2^k-1) 
  \end{align*}
(by counting   in two ways  the pairs $(z,\Delta'')$  with  
 $z\in \Delta''-\Delta^\perp$ and
 $\Delta''\in \Delta^G- (\Delta^\perp  \cap \Delta^G)  $).
Also, $|\Delta|\, \big|\, |\Omega|$ implies that $r\ne \delta\, |\,r$.  But
this condition together with  
$ |\Delta^\perp -\Delta|\ge \{ |\Omega -\Delta^\perp |  /  \! ( 2^\d-2^1 )\}(2^1-1)$
 never hold.
 
Thus, $\Delta'\in \Delta^G$ and $\Delta ' \cap \Delta^\perp\ne 0$ imply that $\Delta'\subseteq\Delta^\perp$,  so   $\Delta^\perp$ is partitioned by 
$   \Delta^\perp  \cap    \Delta^G$.
Finally, if $\Delta'\in \Delta^G$ then either 
$\Delta'\subseteq  \Delta^\perp\subseteq x^\perp$ or 
$\Delta'\cap  \Delta^\perp=0$,
$\<  \Delta', \Delta^\perp\>=V$ and 
$\Delta^\perp\subseteq x^\perp$, so  $\Delta'\not \subseteq x^\perp$. 

\begin{corollary}
If $W $ is an intersection of subspaces $ ( \Delta^g)^\perp,$  $g\in G,$
then $W$ is partitioned by $W\cap \Delta^G$.
\end{corollary}
   
\begin{lemma} 
Either
\begin{itemize}
 \item[\rm (i)]    there is a subspace $\Delta' \in  \Delta ^\perp   \cap   \Delta^G, $
 $\Delta' \ne \Delta,$
 or
  \item[\rm (ii)]    $\dim\Delta=4,$   $V$ has type   $O^+(8,2)$
and  $G\cong A_9 $ is unique up to conjugacy in $\Omega^+(8,2)$.
\end{itemize}
\end{lemma}

  \noindent{\em Description of the   example in}  (11.3ii).
  Let $W=\GF(2)^9$ be the permutation module for $H=A_9$ over~$ GF(2)$,
   and let $\wt(v)$ be the number of nonzero coordinates of  $v\in W$. Then
$V=\{{v\in W\mid }\, \wt(v)\equiv0\pmod{2}\}$ is an $O^+(8,2)$-space with  quadratic form $\phi(v) \equiv \frac{1}{2}\wt(v)\pmod{2}$. 
Clearly, ${H=H'<\Omega^+(8,2)}$,
and   $S=\{\<v\> \mid v \in V, \wt(v)=8\}$ is an $H$-orbit of  $9$ pairwise nonperpendicular points.   
Applying a triality automorphism $\tau$ ([22]; discussed  in 
Section~10~following the proof of (10.2))  produces the desired    $H^\tau$-invariant  set $S^\tau$  of t.s.~4-spaces in  (11.3ii).%

\smallskip \smallskip

{\noindent \em Proof of Lemma.}  Assume that  (i) does not hold.  
  By (11.1), $\Delta^\perp = \Delta $, so $\Delta$ is a maximal t.i.~or t.s.~subspace and   $\d=r$.

 For distinct $\Delta, \Delta'\in \Delta^G$, 
  $G_{\Delta\Delta'}^\Delta$ is antiflag transitive.  (For, if $x\in \Delta$
then $G_{x\Delta'}$ is transitive on the points $y\in 
\Delta'-x^\perp$ and hence on the 
hyperplanes $y^\perp\cap \Delta$ of  $\Delta$ not on $x$, as asserted.)
  Since $\D$ is a minimal block, $G_ \D^\D$ is primitive.  
Moreover, $G_\D$ is  transitive on   
$ \Delta^G-  (  \Delta ^\perp   \cap   \Delta^G )=\Delta^G-\{\D\}$, where 
 $ |\Delta^G|-1=q^{2r-1+c}/|\Delta'-(x^\perp\cap \Delta' )|=q^{2r-1+c}/q^{\d-1}=q^{r+c}$  for $c$ in (9.2).   
Thus,  $G_{\Delta\Delta'}^\Delta$ is primitive by (6.4),  and hence  is as in Theorem I or II.

For  $g\in  G_{\Delta\Delta'} \cap GL(V)$  of order  $p$, let 
 $k=\dim  C_\Delta(g)$, in which~case $g$ also 
 centralizes a $k$-space  of $\Delta'$ 
 (since $\Delta$ and $\Delta'$ are dual $\<g\>$-modules), as well as the anisotropic $(n-2r)$-space 
 $\<\Delta,\Delta'\>^\perp$ unless $V$ has type $O^-(2r+2 ,q)$ with $q$ even, in which case $g$ centralizes at least a 1-space of $\<\Delta,\Delta'\>^\perp$.
 We claim that $k\le (r+1)/2 $.  For otherwise,  if
 $g$ centralizes $\<\Delta,\Delta'\>^\perp$ then 
 $\dim C_V(g) \ge  2k + (n-2r)  > n-r $;  while if $g$ centralizes a 1-space of the 
 anisotropic 2-space $\<\Delta,\Delta'\>^\perp$ then  once again
 $\dim C_V(g) \ge  2k + 1  > n-r $.
Thus,  $C_V(g) $  meets every member of   $ \Delta^G$ nontrivially
 and hence in  a $k$-space.  Now $C_V(g)$ is a subspace 
 having a non-zero t.i.~or t.s.~radical since $|g|=p$, and  having 
 exactly 
 $(q^{r+c}+1)(q^k-1)/(q-1)$ t.i.~or t.s.~points  with  $k >(r+1)/2 $, which is impossible.

Since $\d=r>2$, it follows that     $G_{\Delta\Delta'}^{\,\Delta}$  cannot contain nontrivial transvections;
and it cannot contain $G_2(q)$ by  (5.4d). 
By Theorems I and II,  the only remaining possibilities are  $\delta=4$, $q=2$, and
\vspace{1pt} 
$G_{\Delta\Delta'}^{\,\Delta}=A_6$ or $A_7$. Now $G$ acts on $\Delta^G$ as a
2-transitive group of degree  $q^{r+c}+1= 2^3+1$ or  $ 2^5+1$  
(for $V$ of type    
 $O^+(8,2)$  resp.  $O^-(10,2)$)  in which the stabilizer of two points has a homomorphic image $A_6$ or $A_7$.   In the $O^+(8,2)$ case the action of $A_6$ on $V$   in the preceding {\em Description} readily 
 yields conclusion (i).
In the $O^-(10,2)$ case $W=\< \Delta,\Delta' \>$ is an $O^+(8 ,2)$ space.  
 The action of $A_6$ on $W$   in the preceding {\em Description} shows that there is a third t.s. 4-space of $W$ fixed by $A_6$;
 this must be in $W\cap \Delta^G$, so  $G_{\Delta \Delta' }^W$ cannot be $A_7$
 and hence is $A_6$.  Then $G_{\Delta \Delta' } =A_6$:
  each element of $\Delta^G$ meets $W$ and hence is fixed by $C_G(W)$.
 If $E$ is a Sylow   11-subgroup of $G$, then 
 $|N_{\Omega^-(10,2)}(E)|=33\cdot 5$, so $|N_G(E)|\big| 11\cdot15$, producing the contradiction that 
the number of Sylow 11-subgroups of $G$  is not $ \equiv 1$ (mod 11).

  \smallskip \smallskip 

{\em   From now on we will assume that {\rm (11.3i)} holds.}  Then 
$\Delta $ is not a maximal t.i.~or t.s.~subspace.

  \smallskip \smallskip
  
\begin{lemma}      ~ 
\begin{itemize}
 \item[\rm (i)]   
If $W\ne0,$ $V$  is a   subspace such that
  $W=\<W\cap \Delta^G\>$ and $W^\perp=\<W^\perp\cap \Omega\>,$
 then $W$  is  an intersection of subspaces $ ( \Delta^g)^\perp,$
 $g\in G,$  and is  partitioned by $W\cap \Delta^G$. 
 \item[\rm (ii)]     
If $W_i=\<W_i\cap \Delta^G\>$ for $i=1,2,$  
then $\<W_1,W_2\>\cap \Delta^G$  partitions $\<W_1,W_2\>$.
 \end{itemize}

\end{lemma}

\proof 
(i)  Since $W=\<W\cap \Delta^G\>$, $W^\perp$ is an intersection of subspaces 
 $ ( \Delta^g)^\perp,$
 $g\in G$, and hence by (11.2) is partitioned by subspaces 
$\Delta^g$.  
Since the bilinear form defining the geometry is nondegenerate 
 (cf.~the beginning of Part II),  $W=W^\perp{}^\perp$ is an intersection of subspaces 
$(\Delta^g)^\perp$, and (11.2) applies again.

(ii)    The set $S$ of points  of $ \<W_1,W_2\>$ lying in a member of 
$\<W_1,W_2\>\cap \Delta^G$ clearly spans  $\<W_1,W_2\> $.
We claim that, if $s_1,s_2\in S$ are perpendicular then $\<s_1,s_2\>\subseteq S$,
so $\<S\>=\<W_1,W_2\> $.  Let $s_i\in \Delta_i\in \<W_1,W_2\>\cap \Delta^G$
for $i = 1,2$.  
Since $\d<r$ by (11.3i),   the subspace $ \<\Delta_1,\Delta_2\> ^\perp$ 
is spanned by its points;   by (i),     $\<\Delta_1,\Delta_2\>$
is partitioned by  $\<\Delta_1,\Delta_2\>\cap \Delta^G$, so
$\<s_1,s_2\>\subseteq S$.

\begin{lemma}
$V$ is   orthogonal and $\d = 2$.
 \end{lemma}
\proof
Choose $\Delta'\in \Delta^G,\Delta'\not\subseteq\Delta^\perp.$ 
Then $\Delta'\cap\Delta^\perp=0$ by (11.1), so $W=\<\Delta,\Delta'\> $ is nonsingular.  
By (11.3),  $\delta<r$, so  $W^\perp=\<W^\perp  \cap \Omega\>$  and
$W^\perp \cap  \Delta^G $    partitions $W^\perp$  by   (11.4i).
If $y_1,y_2\in W-x^\perp$, and $y_i\in \Delta_i\in \Delta^G$ ($i=1, 2$), then 
$\Delta_i\subset W$  by (11.4i),
so $W=\<\Delta,\Delta_i\>$ $(i = 1,2)$;
an element of $G_x$ mapping $y_1$ to $y_2$ also maps $\Delta_1$ to $\Delta_2$ and so fixes $W$.
Then $G_W^{\,W}$ is antiflag transitive and imprimitive.
If we are in case (11.3i) for  $G_W^W$  then 
induction implies that $\dim W = 2\d=4$.  The possibility that $W$ is a 4-dimensional symplectic space was excluded at the start of Section 9,
while the 4-dimensional unitary possibility  is eliminated by Kantor-Liebler [14, (5.12)].
Thus,   $V$ is orthogonal. 

It remains to consider the possibility (11.3ii)  for $G_W^{\,W}$,
where we are assuming that $\dim V>8$.     Let 
$\Delta_1\in   W^\perp \cap  \Delta^G$,  
$ \Delta_1'\in (W^\perp \cap  \Delta^G) - \Delta_1'{}^\perp $,
  $W_1=\<\Delta_1, \Delta_1'\>$ and $V'=\< W, W_1\>$,
  so $\dim V'=16$.   
   By  (11.4ii),   $V'$ is  partitioned by $V'\cap \Delta^G$
  Let  $\Delta_2\in V'\cap \Delta^G$, $\Delta_2\not \subseteq W,W'$.
  Then $\Delta_2^\perp \cap  V'$ has dimension 12 and so meets both $W$ and $W'$ in subspaces of dimension $\ge4$.  By (11.4i),  $\Delta_2^\perp\cap  W$ contains some 
  $\Delta_3\in \Delta^G$ and $\Delta_2^\perp\cap  W'$ contains some 
  $\Delta_3'\in \Delta^G$.  Then $\<\Delta_2, \Delta_3, \Delta_3'\>$ is a 
  t.s. subspace of $V'$ of dimension $>8$, which is not possible.
  This   rules out   (11.3ii). 
   
  \smallskip \smallskip
\textsc {Definition}.  Let ${\mathcal L}$  be the  set of all t.s.~subspaces that are intersections of members of $(\Delta^\perp) ^G$.
 By (11.4i), $\Delta^G\subseteq {\mathcal L}$.  Clearly ${\mathcal L}$
 is closed under intersections.

\begin{lemma}
${\mathcal L}$ is the set  of all t.i.~subspaces of a classical geometry
of type  $U(\half n , q^{})$.  
 \end{lemma}
\proof
By (11.2), each member of ${\mathcal L}$ is partitioned by the members of $\Delta^G$  it contains.  If $M$ is a maximal member of ${\mathcal L}$, then $M$ is a maximal t.s.~subspace. 
(For if $x\in M^\perp-M$, then by (11.2) the member of $\Delta^G$ containing $x$ would be in
 $M^\perp$ and, together with $M$, would span  a member of ${\mathcal L}$
   by (11.4i).) 
 
  Assume that  $r=\dim M>4$.  
  Then exactly as in the proof of  (7.2), 
   $M$ is a projective space with $M\cap \Delta^G$ its set of points and 
   $q^\delta+1=q^2+1$ points per line.
   If $\Delta\not\subseteq M=M^\perp$ then $\Delta^\perp \cap M= \<\Delta  , M^\perp\>^\perp$ has dimension $  n-(2+n -r)=r-2$, so $\Delta^\perp \cap M$ is a hyperplane of our new projective space $M$.

 Note that any $N\in {\mathcal L}$ of dimension  $r-\d=r-2$  lies in at least two maximal members of  ${\mathcal L}$: by (11.4i), those maximal members induce a partition of 
 $N^\perp-N$.  If $M$ and $M'$ are $r$-spaces in ${\mathcal L}$  with nonzero intersection, let $N\subset M$ be an  $(r-2)$-space in  ${\mathcal L}$  with  $ M'\cap N \subset M'\cap M$.
 Then  $N^\perp \cap M=M'\cap M$ since 
 $\<N, N^\perp \cap M\>$ is  t.s.~and contains $\<N, M' \cap M\> =M'$.
 If $M''\ne M'$  is an $r$-space in ${\mathcal L}$  containing $N$, then
 $  M''\cap M\subseteq N^\perp\cap M=M'\cap M$; 
 and  $  M''\cap M\subset M'\cap M$, since otherwise $M''\supseteq \<N,M'\cap M\>=M'$.
 Continuing, we find that there exist disjoint $r$-spaces in~${\mathcal L}$.
    
 It follows from Tits [23] that ${\mathcal L}$  is a classical polar space
 since  $r/2>2 $.
 
 Now if $M$ and $M'$ are disjoint maximal subspaces of 
 ${\mathcal L}$  and $\<M, M'\>\ne V$,
 then there is a member of $\Delta^G$ disjoint from   $\<M,M'\>$.
 So $n=\dim V =2r$ or $2r+2$,  where $r$ is even.   If $n=2r + 2$ then $V$ has type $O^-(2r+2,q)$
 and so has  $ ( q^{r+1}+1)(q^r-1) /(q-1)$ points; then 
 $|\Delta^G|=(q^{r+1}+1)(q^r-1) /(q^2-1)$, and ${\mathcal L}$  is of type $U(r+1,q)$.
 Similarly, if $n=2r$, then $V$ has type $O^+(2r,q)$, and the same argument shows 
 ${\mathcal L}$  has type $U(r,q)$.   In all dimensions the results of 
 Tits [23]  show that the embedding of ${\mathcal L}$ in $V(\half n,q ^2)$ is the natural one.
 
 \smallskip
 Next suppose that $r=4$.  
  Then ${\mathcal L}$  is the lattice of points and  lines of a geometry ${\mathcal G}$.  Arguing as above, we find that ${\mathcal G}$  is a generalized quadrangle with  $s=q^2$,
  and
  $t=q$ or $q^3$  according as $V$  has type $O^+(8,q)$ or   $O^-(10,q)$.
 
 If $\Delta ' \not\subseteq  \Delta^\perp$, then $ \<\Delta,\Delta'\>$ has type
 $O^+(4,q)$   and $|\<\Delta,\Delta'\>\cap \Delta^G| =q+1$
 by (11.4i).  For any $\Delta''\in \<\Delta,\Delta'\>\cap \Delta^G, $ $(\Delta'')^\perp 
 \supseteq\Delta^\perp  \cap (\Delta')^\perp$.
  Thus,  if  $t=q^3$ then a theorem of Thas [21] and its proof  identify the quadrangle   as that of type $U(5,q )$,  with uniqueness of the embedding.

  If $t = q$, the points and lines of the quadrangle are certain lines  and
solids of the $O^+(8, q)$  geometry. Any two of the solids are disjoint or meet in a
line,   so they all belong to the same class. Applying the triality map (cf. Section
10), the dual quadrangle is embedded as a set of points and lines in an $O^+(8,q)$
geometry,
  satisfying the hypotheses of Buekenhout-Lef\`evre [1, Theorem 1].
 Thus the dual of ${\mathcal L}$  is of type 
$O^-(6, q)$ in its natural embedding, 
and   ${\mathcal L}$    is of type $U(4, q) $ also embedded naturally.
 This proves (11.6).

 \smallskip 
 
\emph{We can now  complete the proof of} Theorem IV.  By (11.6),  
 ${\mathcal L}$ is    embedded naturally in a projective space derived from a  vector space
 $V(\half n ,q^2)$.  
 Proceeding as in Section 7,  we obtain the original space $V$ by restricting the
scalars;  repeat the argument in that section ({Proof of} (6.1), second paragraph) to show  that   either $q=\d=2$ or  
$e=\delta=2, $  $ q=4$, and that $G$ is
primitive and antiflag transitive on the  $U(\half n ,q)$ geometry.    
If $\half n \ge6$, then this geometry has 
rank $\ge3$; Section 10 does not provide any unitary examples so  (9.4iv) implies that 
$ SU(\half n,q  ) \noreq G\le \Gamma U(\half n,q )< \Gamma O^\epsilon( n,q)$, as required. 
 If $\half n = 4$ or 5 then $G\le \Gamma U(\half n ,q )$ with $q=2$ or 4, and 
 $G \roneq S U(\half n ,q )$ by   Kantor-Liebler [14, (5.12)].   (As in Section 7,
 $G$ must contain $q=2$ or 4 field automorphisms in order to have $G_\Delta^  \Delta\,$   antiflag transitive.)

\section{Rank 4 subgroups of rank 3 groups}
In this section, $G$ will denote a primitive rank 3 permutation group on a set $ X,$
and $H$ a subgroup of  $G$ having rank 4 on $X.$

Let $k,l,\lambda, \mu$  be the usual parameters for $G$, as defined in Higman [9],
and let $I,A,B$ be the adjacency matrices corresponding to the orbits
$\{x\}$, $ \Delta(x)$ and $\Gamma(x)$ of $G_x$, $x\in X$.
If $k,r,s$ are the eigenvalues of $A$, then $\lambda=k+r+rs,$ $
\mu=k+rs,$ $ k(k-\lambda-1)=l\mu$.

\emph{We assume that $H_x$ splits $\Gamma(x)$ into two orbits 
 $\Gamma_1(x)$ and  $\Gamma_2(x)$}, of lengths 
$j,$ $ l - j$ and with adjacency matrices $C, $  $B - C$ respectively. Set
$j t= | \Gamma_1(x)\cap \Delta(y)|$ for $y\in \Gamma_2(x)$.  
Then, with respect to the $\Delta$-graph, the intersection
numbers for  $H$   are as in the following  diagram.
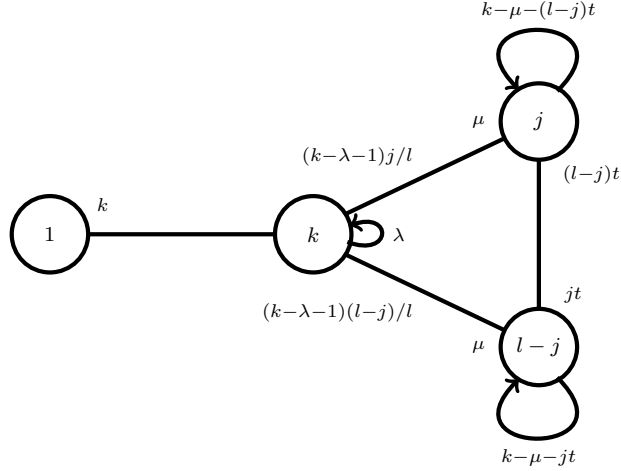
\begin{center}
\begin{tikzpicture}
\draw[ultra thick] (0,0) circle (.2in);
\draw[ultra thick] (-3.5,0) circle (.2in);
\draw[ultra thick] (3,1.5) circle (.2in);
\draw[ultra thick] (3,-1.5) circle (.2in);
\node (0) at (-3.5,0) {\footnotesize$ 1$};
\node (1) at (0,0) { \footnotesize$ k$};
\node (2) at (3,1.5) {\footnotesize$ j $};
\node (3) at (3,-1.5) {\footnotesize$ \ell-j$};
\node (11) at (0.17,0) {$\phantom{X}$};
\node (22) at (3,1.66) {$\phantom{X}$};
\node (33) at (3,-1.66) {$\phantom{X}$};
  \draw[black,ultra thick,-]  (-3.,0) to node [auto] {} (-.5,0);
  \draw[black,ultra thick,-]  (0.41,0.26) to node [auto] 
  {$\scriptstyle (k-\lambda-1)j/\ell$} (2.55,1.28);
  \draw[black,ultra thick,-]  (0.41,-0.26) to node [below left] 
  {$\scriptstyle{(k-\lambda-1)(\ell-j)/\ell}$} (2.55,-1.28);  
 \draw[black,ultra thick,-]  (3,-1) to node [auto] {} (3,1);  

\draw        (11)   [black, ultra thick,out=340,in=20, looseness=30] edge[loop] node [right] {$\scriptstyle \lambda$} (11);
\draw        (22)   [black, ultra thick] edge[loop] node [above] 
{$\scriptstyle k-\mu-(\ell-j)t$} (22);
\draw        (33)   [black, ultra thick,out=315,in=225] edge[loop] node  [below] 
{$\scriptstyle k-\mu-jt$} (33);
\node  (l1) at (-2.8,0.4)  {$\scriptstyle k$};
\node (l3) at (2.2,1.5) {$\scriptstyle \mu$};
\node (l4) at (3.7,0.8) {$\scriptstyle (\ell - j)t$};
\node (l5) at (3.6,-0.8) {$\hspace{-8pt} \scriptstyle jt$};
\node (l6) at (2.2,-1.5) {$\scriptstyle \mu$};
\end{tikzpicture} 
\end{center}

\noindent 
Then $AC=(k-\lambda -1)(j/l)A + (k-\mu - (l-j)t)C +j t(B-C).
$
Applying this to an eigenvector of $A$ and $C$ with eigenvalues $r, \theta, $  respectively,
yields
$$
r\theta =  (k-\lambda -1)(j/l)r+(k-\mu -(l-j)t)\theta +j t(-r-1-\theta).
$$
(Since $A+B+I$ is the all $-1$ matrix, $-r-1$ is an eigenvalue of $B$.)  
Simplifying,
$$
(r(s+1)+lt)\theta = -(j/l)(r+1)(r(s+1)+lt).
$$
Similarly, if $\phi$ is an eigenvalue of $C$ corresponding to the eigenvalue $s$ of $A$, 
$$
(s(r+1)+lt)\phi = -(j/l)(s+1)(s(r+1)+lt).
$$

But the centralizer algebra of $H$ has dimension 4, so exactly one of the eigenspaces
of $A$ must split into two eigenspaces for $C$. If this corresponds to $r$, then  $\theta$ 
is not unique, so
\begin{equation}
r(s+1)+lt = 0.
\end{equation}

\noindent
Since $r\ne s$, it follows that 
\begin{equation}
\phi = -j(s+1)/l.
\end{equation}
But $\phi $ must be an integer, so
\begin{equation}
\mbox{$l/(l , s+1) $  \ divides \ $j$.}
\end{equation}

Also, for $y\in \Gamma_1(x),$  $|\Delta(x)\cap \Gamma_2(y)|= j (l-j)t/k$, so
\begin{equation}
\mbox{$kl$  \ divides \ $j(l-j)r(s+1)$.}
\end{equation}

{\em Remark.}  Of  course, the same results hold in a more general situation (involving association schemes).

\section{Theorem V}
The proof of Theorem V follows (and was inspired by) the pattern of Perin's
Theorem [20] discussed in Section 8. Suppose that $ G$ satisfies the hypotheses of
Theorem V.   If   $G_x$  is transitive on the points outside $x^\perp$, then $G$ is antiflag
transitive. and Theorem IV applies. So we may assume that $G_x$  is transitive on
$x^\perp - x$ and splits $V - x^\perp$  into two orbits.
Then $G$ is transitive on t.i. or t.s. lines.  We use the notation of the last
section.

Suppose first  that $G\le \Gamma Sp(2m,q)$. 
One or both of  $q^{m-1}-1$ and $q^{2(m-1)}-1$ have a primitive divisor $r$
(see (2.4));  let $R\in \Syl_r(G_x)$. Then $W=C_V(R)$ is a nonsingular 2-space and 
$N_G(T)^R$ has rank at most 3.   If $G_x $ has two orbits on the nonsingular
2-spaces containing $x$, then the stabilizer of any projective line (singular or not)
acts 2-transitively on it. By (4.1), $G$  is antiflag transitive, contrary to assumption.
So $G_x $   is transitive on the $q^{2m-2}$  nonsingular 2-spaces containing $x$, and
$G_W^W$ has rank 3;  call the  subdegrees $1,h,q-h$.  As in (8.5),  $(q,h)=(2,1), (3,1),  (4,2), (5,1)
$ or  $(9,3).$

We have $k=q(q^{2m-2}-1)/(q-1),$ $ \ell = q^{2 m-1},$  $ j=q^{2m-2}h$.  Also, 
$r,s=\pm q^{m-1}-1$.  By (12.4),
$$
\mbox {$q^{2m}(q^{2m-2}-1)/(q-1)$ \  divides
 \ $q^{4m-4}h(q-h)q^{ m-1 }(q^{m-1}\pm1),$}
 $$
whence
$$
\mbox {$q^{m-1}\mp1$ \  divides   \ $(q-1)h(q-h).$}$$
This is impossible if $m\ge4$; and none of the specific values of $q$ and $h$ satisfy it when $m=3$.  So this case cannot occur.

 \smallskip
The case $V$ unitary is ruled out by Kantor-Liebler [14, (6.1)].

Suppose  $G\le \Gamma O(2m+1,q), $ $m\ge3$, $q$ odd.  Let $r$ be a primitive divisor of 
$q^{2m-2}-1$, and $R\in \Syl _r ( G_x)$.  Then $W=C_V(R) $ is a nonsingular 3-space, and $N_G(R)^W$  has rank 2 or 3.
(Rank 4 does not occur since $W$   does not contain any t.s. line and hence does not contain any point  of $x^\perp-x$.)
  If $q>5$ then $N_G(R)^W$ contains $\Omega(3,q)$ or (if $q=9$) $A_5$,
  using [5, Chap.~12], and hence so does $C_G(W^\perp)^W$
  since  $N_G(R)^{W^\perp}$ is solvable.    Then the argument used  in (6.2) (i.e., using a group behaving like $Q$)    shows that
 $G\ge \Omega(2m+ 1,q)$, 
  which is a contradiction since that group has rank 3 on points.
  Thus, $q\le5$.
 \smallskip

Suppose  $G\le \Gamma O(2m+1,5), $ $m\ge3$.  
In addition to $r$ we  will use  a primitive divisor 
 $r^\circ$
 of $5^{m-1}-1$; let  $R^\circ\in \Syl_{r^\circ}(G_x)$   and $W^\circ=C_V(R^\circ)$.
As above, we may assume that $N_G(R)^{W}$
and $N_G(R^\circ)^{W^\circ}$ 
do  not contain $\Omega(3,5)$;
both are   rank~3  groups that therefore contain    $S_4$.
Since 
    $N_G(R)^{W^\perp} \le \Gamma O^-(2, 5^{m-1})$ 
    and 
    $N_G(R^\circ{})^{W^\circ{}^\perp} \le \Gamma O^+(2, 5^{m-1})$
   are  metacyclic,  $C_G(W^\perp)\cap \Omega(V)$ and  
    $C_G(W^\circ{}^\perp)\cap \Omega(V)$  
   contain normal subgroups $A$ and $B$, respectively,  
   isomorphic to $\Z_2^2$.  In view of the behavior of $R^{W^\perp}$ and 
$R^\circ{}^{W^\circ{}^\perp}$, $W$ and $W^\circ$ are not isometric.   
   
   Let $b_1,b_2,b_3$ be an orthogonal basis of $W$ with respect to which $A$ is diagonal;  $N_G(R)$ acts transitively on $\{\<b_1\>,\<b_2\>,\<b_3\>\}$.  
   Then  $T=\<b_1+2b_2,b_2+2b_3\>$ is an $O^+(2,5)$-space.
   Since $N_G(R^\circ{})^{W^\circ{}^\perp} $ contains representatives of both $G$-orbits of $O^+(2,5)$-spaces, we may assume that 
   $T\subset  W^\circ   $.   Then $F=\<W, W^\circ \>$ is a 4-space 
   containing  non-isometric nonsingular 3-spaces, and hence is nonsingular. 
   We have a  group $H=\<A,B\>^F\le\Omega(F)$;
   $A$ and $B$ are not conjugate (since $[V,A]=W$, 
    $[V,B]=W^\circ$), and   $A\cap B=1$ (every nontrivial element of $A$ moves $T$).
     However, $\Omega(F)=\Omega^\pm(4,5)$ has no such subgroup  $H$.

   \smallskip

The case  $G\le \Gamma O(2m+1,3), $ $m\ge3$, is   harder. 
This time choose a 
primitive divisor  ${r\mid 3^{2m-2}-1}$ or 
 $r\mid 3^{ m-1}-1$ according to whether  $m$ is  odd or even.
  Then $W=C_V(R) $ is a nonsingular 3-space;  $N_G(R)^W$ 
   has rank 2 or 3 and so contains $\Omega(3,3)$ or $D_8$, while 
    $N_G(R)^{W^\perp} \le \Gamma O^\pm(2,3^{m-1})$ is metacyclic, with a normal cyclic subgroup of order dividing $3^{ m-1}\pm1\equiv2$  (mod 4)
    in view of our choice of $r$.
If $t$ is the square of an element of  order 4 in $N_G(R)$, 
or  if $N_G(R)^W \ge \Omega(3,3)$ and $t\in N_G(R)'$ has order 2,  then
  $t\in C_G(W^\perp)$  is an involution with $t^W $
 inducing $-1$ on an anisotropic  2-space.
 Let $b,b' \in W$ be linearly independent vectors with $b^t=-b$, $b'^t=-b'$,
 where   (since $q=3$ and in view of   the action $N_G(R)^W$) we may  assume that they are perpendicular
 and ${\phi(b)=\phi(b')=1} $  for the quadratic form $\phi$ on~$V$.

  Any two   $G_{\<b\> }$-conjugates of  the reflection $-t^{b^\perp}$ commute.
 (For otherwise, the product of two such non-commuting conjugates of $t$ has
order 3 and centralizes   $y^\perp/y$ for some point $y\in b^\perp$.  Then the  argument
used  in (6.2)     yields the contradiction  $G\ge \Omega(2m+1,3)$.)
It follows that 
 $b^\perp=U_1\perp U_2$, where $U_1=\< b'^{G_{\<b\>}}\> $
 is spanned by pairwise perpendicular members of 
 $\< b'\>^{G_{\<b\>}}$.
  
 Let ${\mathcal N}_1
  =\{\<v\>\mid  v\in V, \, \phi(v)=1\}$.
 Since $G_x$ has two orbits  of   $y\in \Omega -x^\perp$
 and each $\<x,y\> $ contains a unique   $\<b\> \in 
 {\mathcal N}_1$,
  $G_x$  has two orbits  on $  {\mathcal N}_1- x^\perp$.  This proves that
   $G_{\<b\>} $
 has at most two orbits on  $ \Omega - b^\perp;$
 and   there are two orbits if and only if  $G$ is transitive on  
 $  {\mathcal N}_1$.
 
Suppose that  $G$ is intransitive on $  {\mathcal N}_1$.    Then  $G_{\<b\>}$ is transitive
on $\Omega-b^\perp$, but leaves invariant $U_1$ and $U_2$.  Then $U_2=0$ and 
 $G_{\<b\>}$ is monomial on  $U_1= b^\perp$ with respect to an orthonormal basis.
 Since $G_{\<b\>}$ is transitive
on $\Omega-b^\perp$ and $2m=n-1>4$, this is impossible.

Thus, $G$ is transitive on $  {\mathcal N}_1$.  
 Let $s$ be a primitive divisor of $3^m-1$ or $3^{2m}-1$
 such that  $s\big | |  \Omega\cap  b^\perp |$.
 If $S\in \Syl_s(G)$ then $C_V(S)\in {\mathcal N}_{1}$.
 We may assume that $S$ fixes $\<b\>$ and hence has  no proper nonsingular invariant subspace
 $U_2$ in  $ b^\perp$.    
 Once again $G_{\<b\>}$ is monomial on  $U_1=b^\perp$ with respect to an orthonormal basis.
  Members of  $ \Omega -b^\perp $  look like $\<b+u\>$ with 
  $u\in b^\perp $  and $\phi(u)=-1$, where   $u$ has $k$ nonzero coordinates
 with $k\equiv 2$ (mod 3).  Since there are only two such orbits, $k$ can only be 2 or 5, so  $\dim b^\perp=n-1<8$ and   we are in an $O(7,3)$ geometry.
Since   
$G$ is transitive on   ${\mathcal N}_1$,   $G_x$ has  an orbit  on 
${\mathcal N}_1 - x^\perp$    of length
 $
 \mbox{$   \frac{1}{2}$} 3^{3}(3^3\pm1)\cdot
 2^2\binom{7-1}{2}\big/ \!\{(3^{6}-1)/(3-1)\}\!,  
$
 which   is not an integer.%
\smallskip 

     Finally, consider the case 
    $G\le \Gamma O^\pm(2m,q), $  $m\ge3$, where $q>2$ (by hypothesis),  in which   $x^\perp/x$ has  
$(q^{m-1}\mp1)(q^{m-2}\pm1)/(q-1)$ points.  If $m=3$, use Kantor-Liebler [14, (5.12) and (5.14)].
Assume that $m\ge4$, and use $r\mid q^{m-2}\pm1$ and $R\in \Syl_r(G_x)$ as before,
temporarily excluding the case $O^-(8,q)$ with $q$ a Mersenne prime.
This time  $W=C_V(R)$ is a nonsingular 4-space
of type $O^-(4,q)$ since $V$ has type  $O^\pm(2m,q)$  and $[V,R]$  has type $O^\mp(2m-4,q)$.
Then $N_G(R)^W$ has rank 2 or 3 and hence  contains $\Omega^-(4,q)$ or (if $q=3)\, $ $A_5$ [5, Chap.~12].  As in  (6.2) we obtain the contradiction  $G\ge \Omega^\pm(2m,q)$.
  
  This leaves the  excluded   possibility 
  $G\le \Gamma O^-(8,q)$ with $q\,$ a Mersenne prime.  
  We may assume that $-1\in G$.  
 If  $L$ is a line then $G_L^L$ is 2-transitive and hence contains $SL(2,q)$. 
   Then there is an involution $t\in G$ such that $-t=1$ on $L$ and 
   $W=C_V(-t)$   has type $O^+(4,q)$.  Let $-t\in R\in \Syl_2(C_G(L))$,  so $W=C_V(R)$ and  $R\in \Syl_2(C_G(W))$.
  By the Frattini argument, if $N=N_G(R)$  then $N ^W$   is transitive on lines while 
  $N_L^L$ is 2-transitive.
  Then  $N ^W\ge  \Omega^+(4,q).2; $
  clearly $N^{W^\perp}\!  \le \Gamma O^-(4,q)$.
Thus,  if $q>3$ then $C_{N }(W^\perp)$  contains $ \Omega^+(4,q)$, 
hence a long root group, 
and then all long root groups by line-transitivity; but this produces the usual 
contradiction $G\ge \Omega(V)$.  If $q=3$ then $C_{N }(W^\perp)$ contains
 an involution centralizing a  6-space,  and   a simpler version of the argument used 
above for   $\Omega(2m+1,3)$ produces a contradiction.
 This   completes the proof of Theorem V.

\smallskip   \smallskip
{\em Remark.}  If $G < O^\pm(2m, 2$), the argument breaks down when
 $r \mid 2^ { m- 2} \pm1$,  $\dim W = 4$, and  $|N_G(R)^W | = 10$ or 20.

\section{Concluding remarks}
1. The method used in our proofs for employing $p$-groups also works for
suitable permutation representations of the exceptional Chevalley groups.
\smallskip

2. After classifying antiflag transitive groups, it is natural to ask  about
transitivity on incident point-hyperplane pairs (where the hyperplane is not the
polar of the point in the case of a classical geometry). If a group $G$  is transitive
on all such  pairs in $PG(n - 1, q)$, then it is transitive on  incident point-line
pairs, and hence 2-transitive on points (Kantor [12]); so Theorem I applies.
However, for classical geometries, results are known only  in the  unitary case
(Kantor-Liebler [14, (6.1)]).%
\smallskip

3.
The proofs of Theorems I-III do not depend  on ``modern'' group-theoretic classification theorems.
 Theorem~IV requires results summarized  in Kantor-Liebler [14] that only  use older  group theory; most of  
the required results used nothing more than  elementary  properties of classical groups, such as concrete sets of generators.%
 
\smallskip
4. It should be noted that [14] produces a proof of the rank 2
analogue of Theorem IV, as follows.   We assume that $V$  does not have type
$O^+(4, q)$. The primitive case proceeds as in Sections 9, 10. In the imprimitive
case, the block $\Delta$  of Section 11 is a t.i.~or t.s.~line. If $x\in \Delta $ then $G_x$ is transitive
on $\Delta^G - \{\Delta\}$   and hence on   $x^\perp -\Delta$. It follows easily that
$G$   has one orbit of points and two
orbits each of lines and incident point-line pairs.   Now [14, Sect. 5] applies.

\appendix
\section{The $G_2(q)$ generalized  hexagon}

This appendix contains new and elementary proofs of the existence and
uniqueness statements in Section 3, as well as further properties of the generalized 
hexagons (including antiflag transitivity).

Assume that ${\mathcal G}$ is as in (3.2), and set $W(x) = W_1(x)$.  We will prove several
properties of ${\mathcal G}$, from which an explicit construction will easily follow.

\begin{lemma}
{\rm (i)}    For any points $x,y $ of ${\mathcal G}$ such that 
$d(x,y)=1$ or $2,$   all $1$-spaces of 
$\<x,y\>$  are points of  ${\mathcal G}.$

\begin{itemize}
\item[\rm (ii)]    If $z\notin \<x,y\>$ and $d(x,y)=1$ or $2,$  then $W_2(z)\cap \<x,y\>$ is either 
$\<x,y\>$  or a point.
\item[\rm (iii)]   
 Either $\dim V=6$ and $V$ is symplectic$,$ or $\dim V=7$ and $V$ is orthogonal.
In either case$,$ the  points and lines of ${\mathcal G}$ consist of all points and certain t.i.~or t.s.~lines of $V.$  Moreover$,$ $W_2(x)$ consists of all points of $x^\perp$ $($i.e.$,$ $d(x,y)\le2\iff  y\in x^\perp)$.
\end{itemize}
\end{lemma}
\proof
 (i)   $\<x,y\>\subseteq W(u)$ if $u\in W(x)\cap    W(y)$. 
 \smallskip
 
 (ii) If $d(x, y) = 1$, this follows from the axioms for a generalized hexagon.
Suppose  $d(x, y) = 2$, and set $u= W (x)\cap W (y)$. 
We must show that the subspace $W_2(z)\cap \<x,y\>$ is nonzero (cf. (f) in Section 2).
This  is clear if $d(u,z)\le2$, while if $d(u,z) =3$  it follows from the fact that
$ W(u) \cap W_2(z)$ is a subspace meeting each line on $u$. 
\smallskip
  
 (iii)    As in Yanushka [27, Sect. 3], this follows   from (ii): the points and lines are
 the points and lines of a polar space (Tits [23]).
 Moreover,  $\mathcal G$    has exactly 
  $(q^6 - 1)/(q - 1) $  points  (and  $|W_2(x)|=(q^5-1)/(q-1)$).
  
  \smallskip
 \smallskip
 
 Two points are {\em opposite} if they are at distance 3.

\begin{lemma}
Let $a$ and $b$ be opposite points$,$ and set $H = \<W(a), W(b) \>$.
\begin{itemize}
\item[\rm (i)]    $H = E \oplus F ,$ where $E$ and $F$ are t.i.~or t.s.~planes such that$,$ for $e \in  E,$ 
$f\in F, $  $ \<e,f\>$ is a ${\mathcal G}$-line if and only if it is a $($t.i.~or t.s.$)$ line 
$($call these $E|F$-{\rm lines}$).$
\item[\rm (ii)] 
If $e \in E,$ then $W(e) = \<e, e^\perp \cap F\>.$
\item[\rm (iii)]  If $x $  is a point on no $E| F$-line$,$ then $W(x)$ meets exactly $q +1$
$ E | F$-lines$,$
and the points of intersection lie on   a t.i.~or t.s.~line.
\item[\rm (iv)] 
If $V$ has type $Sp(6, q),$  then $q$ is even.

\end{itemize}
\end{lemma}

\proof
(i) Since $W(a) \cap W(b) = 0$,   $\dim H = 6$.  
Let $a = x_1 , x_2 , x_3 ,
b = x_4, x_5,x_6$   be the vertices of an ordinary hexagon in ${\mathcal G}$. Then
$x_2,x_6\in W(a)$ and $x_3,x_5\in W(b)$. 
Set 
$E=\<x_2,x_4,x_6\>$ and $F=\<x_1,x_3,x_5\>$.~Then 
$E$ and $F$ are t.i.~or t.s.~(by (A.1iii)) and $H=E\oplus F$.
Also $W(x_{2i}) =\<x_{2i}, x_{2i}^\perp \cap F\>$    for each $i$.
We can thus
vary $x_2 , x_6\in  W(a) \cap   E$, and also move around the ordinary hexagon, in order
to show that each t.i.~or t.s.~line $\<e,f\>$ is a ${\mathcal G}$-line
(for $e\in E, f\in F$).
\smallskip

(ii) This is clear from the above proof. (In fact, the points of $E \cup  F$ and
the $E| F$-lines form a degenerate subhexagon with $s = 1, t = q.)$
 \smallskip
 
(iii)    
Let $E_x=E\cap x^\perp, F_x=E\cap x^\perp$,  
$U=\<E_x,F_x\>$, $e=F_x^\perp\cap E$ and $
f=E_x^\perp\cap F$.
The pair $e,f$ corresponds to a flag of $E$ (and of $F$) if and only if $e$ and $f$ are perpendicular; and then $e\in E_x, f\in F_x$ and 
(for $V$ symplectic resp.  orthogonal)
$U^\perp $ is 
$\<e,f\>$
or $\<e,f\>\perp H^\perp$, which   cannot contain  the point $x$ lying in no $(E|F)$-line.  

Thus, $e,f$ corresponds to an antiflag of $E$. It follows easily that $U$ is nonsingular.
If $z \in E_ x$ and  $u=W(x)\cap W(z)$ (cf. (A.1iii)), 
then  $ \<z ,u\>$  is a ${\mathcal G}$-line and hence 
(by (ii))
an $(E|F)$-line, so  $u \in W(x)\cap U$. 
Thus, $W(x)\cap U$ is the desired set of points, and is   a t.i. or t.s.  line.

\smallskip
(iv) If $V$ has type $Sp(6, q)$, then $U$ has type $Sp(4, q)$. But the $Sp(4, q)$
quadrangle contains six lines forming a $3 \times 3$ grid (such as $E \cap  x^\perp, $ $ F \cap  x^\perp,$ $ W(x)\cap U$, and any three $E |F$-lines in $U$) if and only if $q$ is even.

\smallskip\smallskip
{\em Remark.}  Because of (A.2iv), and the isomorphism between the $Sp(6, q)$
and $O(7, q)$ geometries when $q$ is even, we will assume from now on that $V$ has
type  $O(7, q)$.  Then $H$ has type $O^+(6, q)$, and the line mentioned in (iii) is
$W(x) \cap  H$. 
Also, $O(7, q) = SO(7, q) \times \{\pm1\}$, so we may where necessary
assume that linear automorphisms of ${\mathcal G}$ have determinant 1.

\smallskip\smallskip
The next lemma is more technical, and concerns generating ${\mathcal G}$.

\begin{lemma}
Let $S$ be a set of points$,$ containing at least one pair $a, b$ of
opposite points$,$ and such that $W(a) \cap b^\perp \subseteq S $ for any such pair. Then either
$S = E \cup  F$ for some  $E, F$ as in {\rm(A.2i),} or $S$ consists of all points of ${\mathcal G}$.

\end{lemma}

\proof
Certainly  $S\supseteq  E\cup F$  if $E$ and $F$ are obtained as in (A.2i).  
(Each line of $F$ on $b$ is $W(a')\cap F$ for some $a'\in W(b)\cap  E=W(b)\cap  a^\perp$.)
 Let ${\mathcal G}_0$  consist of  $S$ together with the set of   lines meeting it at least twice. We will show that ${\mathcal G}_0$ is a (possibly degenerate)
subhexagon.

Let  $L $  be a line  of  ${\mathcal G}_0$  and $x\in  S -L$; we must show that the unique point $u$ of $L$
nearest $x$ lies in $S.$  Let $y \in L - u$. Since $x$ is opposite some point of $E$
or $F,$ our
hypothesis implies that each line on $x$ meets   $S - \{x\}$. If $d(x, u) = 1$, pick
$z \in  S \cap W(x)$ with  $d(y, z) = 3$, so $u\in W(y)\cap z^\perp\subset  S$.
  If $d(x,u)=2$ then $d(x,y)=3$, so    $u\in W(y)\cap x^\perp \subset S$.

Thus,  ${\mathcal G}_0$  is a subhexagon.  Let $a\in S$.  Then $S\cap W(a)$ has the following properties:   it meets every line on $a$ at least twice; 
if $x,y\in S\cap W(a)$ and $W(a)=\<a,x,y\>$, then $\<x,y\>\subseteq S$.  (For,
 since ${\mathcal G}_0$  is a subhexagon,  there is a point $b\in x^\perp \cap y^\perp\cap S$   opposite $a$, and then 
$\<x,y\> =W(a) \cap b^\perp.$)  Thus $S\cap W(a)$ is a 
subplane of $W(a)$ (possibly degenerate:  just $\{a\}\cup \<x,y\>$).

If each line of ${\mathcal G}_0$  has size 2, then $S=E\cup F$.  So suppose that some line of 
${\mathcal G}_0$ on $a$ has at least three  points.  Then $S\cap W(a)$ is  nondegenerate, and hence is all of $W(a)$.   Thus ${\mathcal G} = {\mathcal G}_0$.

\begin{lemma}
Suppose ${\mathcal G}$    and  ${\mathcal G}'$ are both embedded in $V$ as in Section $3.$ Let
$x_1,\dots ,x_6$ and $y_1,\dots ,y_6$ be the vertices of ordinary hexagons in 
${\mathcal G}$    resp.   ${\mathcal G}'.$  Then there is an element of $GL(V)$ mapping
 $x_i$ to $y_i$  $(i = 1,\dots,6)$ and inducing an isomorphism of 
${\mathcal G}$    onto  ${\mathcal G}'$.
\end{lemma}

\proof
The orthogonal geometries determined by ${\mathcal G}$    and  ${\mathcal G}'$
as in (A.1iii) are equivalent under $GL(V)$; so we may suppose that they are equal.
There is an orthogonal transformation taking $x_i$ to $y_i$
$(i = 1,\dots,6)$, so we may assume that $x_i=y_i$ for each $i$.  Set $E=\<x_2,x_4,x_6\>$,
$F=\<x_1,x_3,x_5\>$. 
 By (A.2ii),  if  $e\in E, f\in F$, then
$W(e)$ and $W(f)$ are the same whether computed in ${\mathcal G}$    or  ${\mathcal G}'$.

Pick a point $x$ on no $E | F$-line,   
so  $W(x)\cap H$ is  the t.s.~line  in  (A.2iii),
and hence is   one of the $q-1$ lines 
$\ne E\cap x^\perp,  F\cap x^\perp$
in  $U$
 meeting each $E|F$-line of 
 $U=\<E\cap x^\perp, F\cap x^\perp\>$.
But $O(7,q)_{EFU}$ is transitive  on these $q-1$ lines, so we may assume that $W(x)=\<x, W(x)\cap H \>$   is the same in 
${\mathcal G}$    and  ${\mathcal G}'$ for the chosen $x$.

We will show that (A.3) applies to 
the set  $S$  of points $u$ of $V$ such that $W(u)$ is the same in both
${\mathcal G}$    and  ${\mathcal G}'$.
 Let $a,b\in S$ be opposite.  
Then $A=W(a)\cap b^\perp$ and $B=W(b)\cap a^\perp$ are t.s.~lines.  
If  $u\in A$ then $L =W(b)\cap u^\perp$ is a line on $b$;  let  $v=L \cap B$.  In 
${\mathcal G}$    (and  ${\mathcal G}'$)  there is a unique shortest path 
$a,u,w,b$;  since $w\in W(b)\cap a^\perp$ and $w\in u^\perp$ we have $w=v$.  Then $W(u)=\<u,a,v\>$ in both
${\mathcal G}$    and  ${\mathcal G}'$, so $u\in S$.

  Now   ${\mathcal G} = {\mathcal G}'$  by (A.3).

\begin{corollary}
The group $\Aut_V({\mathcal G})$ of automorphisms of ${\mathcal G}$ induced by elements of $SL(V)$ is transitive on the set of ordered ordinary hexagons of ${\mathcal G}$.
In particular$,$ $\Aut_V({\mathcal G})$ is antiflag transitive.

\end{corollary}

\begin{corollary}
{\rm (i)}
 There is a subgroup $K\cong SL(3,q)$ of $\Aut_V({\mathcal G})$ fixing $E$ and $F$ $($cf. $(${\rm A.2}$))$ and centralizing $H^\perp$.
\begin{itemize}
\item[\rm (ii)]   The stabilizer of $E$ in $\Aut_V({\mathcal G})$ induces $SL(3,q)$ on it.
 \item[\rm (iii)]      $|\Aut_V({\mathcal G})|=(q^6-1)q^6(q^2-1)$ 
 and $\Aut_V({\mathcal G})$ contains no nontrivial scalar transformations.    
   \end{itemize}
\end{corollary}

\proof
(i) Use  (A.4)  and (2.1)  (compare (5.4c)).%
\smallskip

(ii) 
The plane $E$ uniquely determines the plane 
$F=\<W(a)\cap W(b)\mid a,b\in E,$ $ \,a\ne b \>$.  
Let  $J=\Aut_V({\mathcal G})_{EFU}$ and $C=C_J(E\cap U)$ 
for  the $ O^+(4,q)$-space   $U=\<E\cap x^\perp, F\cap x^\perp\>$   in the proof of (A.2iii)  and (A.4).  
Both  $J$ and $K_U$  fix the antiflag $(E\cap U^\perp,E\cap U)$ of  $E$
and induce $GL(2,q)$ on $E\cap U$.  Then  $J=C  K_U$.
We will show that $C=1$, so that 
$\Aut_V({\mathcal G})_E=\Aut_V({\mathcal G})_{EF} =K$ and
(ii) holds.

Since 
 $C $ is 1 on $E\cap U$, fixes $F\cap U$ and acts inside $O(U)=O^+(4,q)$,
 it is 1 on $U$.  
 Then  $C$ fixes  each 2-space $W(y)\cap U$ for $y\in  U^\perp-\<E,F\>$,
 and then  fixes the unique point $y$ joined     by ${\mathcal G}$-lines
to all points of  $ W(y)\cap U$.  Since $C$ fixes  $U^\perp\cap E$ and $U^\perp\cap F$,
$C<SL(V)$ centralizes $U$  and fixes all points of $U^\perp$, so 
  $C=1$. 
        
   \smallskip
 (iii) 
 There are  $[(q^6-1)/(q-1)]\cdot (q+1)q\cdot qq\cdot  qq\cdot q$ ordered hexagons in ${\mathcal G}$.  The stabilizer in 
$ \Aut_V({\mathcal G})$ of one of them is the stabilizer in $K$ of a triangle in $E$ and hence has order  $ (q-1)^2$. 
 The final assertion is clear since    $\Aut_V({\mathcal G})<O(V)\cap SL(V)$.
 
\begin{theorem}
 Each $O(7, q)$ space has one and only one isomorphism type  of generalized
hexagons embedded as in Section $3$. An $Sp(6, q)$ space has such a hexagon if and
only if $q$ is even.
\end{theorem}

\proof
Uniqueness follows from (A.4), and the assertion about $Sp(6, q)$
from (A.2iv). The preceding results (especially (A.1), (A.2) and (A.6)) tell us
exactly how ${\mathcal G}$  must look, and hence how to construct ${\mathcal G}$.
 \smallskip  \smallskip
 
 {\em Construction.}  Let $V$ be a vector space carrying a geometry of type 
 $O(7, q),$
and $E$ and  $F$   t.s.~planes such that  $H=\<E,F\>$ is nonsingular of dimension $6$.  Let $K<O(7,q)$ fix $E$ and  $F$, centralize $H^\perp$, and induce $SL(3,q)$  on both $E$ and $F$.
 If $\{e_1,e_2,e_3\}$ is a basis for $E$ and $\{f_1,f_2,f_3\}$ the dual  basis for $F$,
 then the matrices of $g^E$ and $g^F$ with respect to these bases are inverse transposes of one another for all $g\in K$.  
We may assume that  $H^\perp=\<d\> $  with $\phi(d) =-1$.

We must use the $E|F$-lines $\<e,f\>$, with $e\in E, f\in e^\perp\cap  F$, as ${\mathcal G}$-lines;
set $W(e)= \<e,e^\perp \cap F\>  , $ $W(f)= \<f, f^\perp \cap E \>$
as in (A.2ii).
Note that $K$  is transitive on the 
$(q^2+q+1)(q+1)(q-1)$ points of $H$  not in $E\cup F$, 
on the $(q^2+q+1)(q^3-q^2)$ points of $V-H$, and
on the 
$(q^2+q+1)(q^3-q^2)$  lines of $H$    not 
meeting $E\cup F$.

We will use the  $E|F$-line $\<e_1, f_2\>$, the point $u=\<e_1+ f_2\> $,
and the t.s.~plane  
$W(u)=\<e_1,f_2,e_3 + f_3 + d  \>$.
Write $W(u^g)=W(u)^g$ for all $g\in K$.
\emph{The new points must be the t.s.~points of $V-H,$} and 
\emph{the new ${\mathcal G}$-lines must be the lines of $W(u^g)$ through $u^g,$ for all $g\in K$.}
 We must show that this is well-defined and yields a generalized hexagon.  This will 
be done in several steps.
\smallskip

(1)  \emph{If $u^g=u$ then $W(u^g)=W(u)$}; so $W(u^g)$ is well-defined.
For, $|K_u|=q^3(q-1)$, and $K_u$ fixes $W(\<e_1\>)/\<e_1, f_2\>$ and $W(\<f_2\>)/\<e_1, f_2\>$.
Thus, each $p$-element of $K_u$ fixes every plane containing $\<e_1,f_2\>$.  Suppose  $|g|\big| q-1$.  
Since $g^E$ and $g^F$ are diagonalizable, we may assume that our dual bases 
$\{e_1,e_2,e_3\}$ and $\{f_1,f_2,f_3\}$
have been chosen so that  $g$ fixes each $\<e_i\>$, $\<f_i\>$.
If $g^E=\diag(\a,\b,\g)$ then $g^F=\diag(\a^{-1},\b^{-1},\g^{-1})$ and
$\a\b\g=1$.   Since $u^g=u=\<e_1+ f_2\>$ we have $\b^{-1}=\a$, whence $e_3^g=e_3,$ $f_3^g=f_3$.
Then $W(u)^g=\<\a e_1,\b^{-1} f_2,e_3+f_3+d\>=W(u)$.

\smallskip  \smallskip

(2)   \emph{If $\,W(u) ^g=W(u)$ then $u^g=u$}.  For,  $g$ fixes 
$W(u)\cap E=\<e_1\>$ and $W(u)\cap F=\<f_2\>$.  Here,  $K_{\<e_1\>\<f_2\>}$ is the stabilizer of a flag of $PG(2,q)$,
 of order  $q^3(q-1)^2$; each of its $p$-elements fixes $u$.  
 If $|g|\big| q-1$  then $g$ is diagonalizable and
 we may assume that our dual bases 
$\{e_1,e_2,e_3\}$ and $\{f_1,f_2,f_3\}$
have been chosen so that  $g$ fixes each $\<e_i\>$, $\<f_i\>$.
Since $g$ fixes  $W(u)=\<e_1,f_2,e_3 + f_3 + d  \>$,
if $g^E=\diag(\a,\b,\g)$~with  $\a\b\g=1$ then 
$ ( e_3 + f_3 + d)^g = \g e_3+\g^{-1}f_3+d$, so $\g=1$, whence
$\b^{-1}=\a$ and~$u^g=u$.

\smallskip   
(3)  \emph{If $L$ is a ${\mathcal G}$-line on $u$ then $L\subset W(u)$}.  (For, we may assume  
$L\not\subseteq H$ and $u^g\in L \subset W(u^g)$ for some $g\in G$, so $u=L\cap H=u^g$ and $L\subset W(u)$.)   The total number of  ${\mathcal G}$-lines is then 
$$
(q^2+q+1)(q+1)+(q^2+q+1)(q+1)\cdot (q-1)q = (q^6-1)/(q-1).
$$
Since $K$ is transitive on $V-H$, each point $x\notin H$ lies on 
$$
(q^2+q+1)(q+1) (q-1) q\cdot q /(q^2+q+1)(q^3-q^2)=q+1
$$
${\mathcal G}$-lines.

\smallskip  
(4)   Let $x\in V-H$.  Then $K_x \cong SL(2,q)$ acts on the 
$O^+(4,q)$-space  
 $U=\<E\cap x^\perp ,F\cap x^\perp \>$;  it  
fixes each of the $q-1$ lines $M\ne E\cap x^\perp ,F\cap x^\perp$ of the same type as $E\cap x^\perp $ that partition the points of 
 $U$, and $K_x^M \cong SL(2,q)$.

If $L$ is a ${\mathcal G}$-line  on $x$ then $y=L\cap H$ is singular but not  in $E\cup F$, and $W(y)$ 
contains $x$ and points  $e\in E$ and $f\in F$.  Since $W(y)$ is t.s. it follows that 
$y\in \<e,f\> $  lies in   $\<E\cap x^\perp ,F\cap x^\perp \>=U$ and hence on one of the lines $M$.  

 \emph{Define}  $W(x)=\<x,M\>$; this is a t.s.~plane.
Since $K_x^M$ is transitive, all  lines of  $W(x)$ on $x$ are 
${\mathcal G}$-lines.
By the transitivity of $K$ on
 the t.s.~lines of $ H$ not  meeting  $E \cup F$, each such line 
 occurs as $W(z) \cap H$   for some
$z\in V - H$.  Since the numbers of such $z$ and such t.s.~lines are the same, distinct
points $z$  yield distinct $W(z).$   It follows that   $W(a)\ne W(b)$ 
\emph{for any   distinct points  $a,b$ of $V.$}
\smallskip 

(5) \emph{Points $a,b$ are perpendicular if and only if $d(a,b)\le 2$}.  For, 
if   $1\le d(a,b)\le 2$ then  $a,b\in W(c)$ for some $c$, and $W(c)$ is t.s.   But the number of such  ordered  pairs is 
$\{(q^6-1)/(q-1)\} \! (q+1)q   +   \{(q^6-1)/(q-1)\} \! (q^2+q)q^2$,
which is  the same as the number of ordered pairs of distinct perpendicular points.

\smallskip  
(6) \emph{${\mathcal G}$ has no $k$-gons for $k\le 5$}.  For, let $a_1,\dots , a_k$ be the vertices of a $k$-gon.
Then $d(a_i,a_j)\le 2$ for all $i,j$ so $\<a_1,\dots, a_k\>$ is a t.s.~plane
by (5), which must be both  $W(a_1)$ and $W(a_2)$, contradicting (4).

\smallskip
(7) \emph{${\mathcal G}$  is a generalized hexagon}.  
Since each ${\mathcal G}$-line is on $q+1$ points, and each point is on 
$q+1\,$  ${\mathcal G}$-lines, 
this follows from the same type of elementary counting argument as in the proof of (10.2).
\smallskip  

This completes the proof of (A.7).
\smallskip\smallskip

{\em Remarks.}   Further properties of the group
$G_2(q)=\Aut_V({\mathcal G})$ are found in (5.4).  Additional information, such as 
simplicity when $q\ne 2$ and identification with $PSU(3,3)\semi Z_2$ if $q=2$, 
are left to the reader, and can be found in Tits [22].

\end{document}